\input amstex
\loadeufm

\documentstyle{amsppt}

\magnification=\magstep1

\baselineskip=20pt
\parskip=5.5pt
\hsize=6.5truein
\vsize=9truein
\NoBlackBoxes

\define\br{{\Bbb R}}

\define\e{{\varepsilon}}
\define\OO{{\Omega}}

\define\bu{{\bold{u}}}
\define\bbf{{\bold{f}}}
\define\bo{{\bold{0}}}
\define\bv{{\bold{v}}}
\define\bg{{\bold{g}}}
\define\bw{{\bold{w}}}
\topmatter

\title
The $L^p$ Boundary Value Problems on Lipschitz Domains
\endtitle

\author Zhongwei Shen
\endauthor

\leftheadtext{Zhongwei Shen}
\rightheadtext{The $L^p$ Invertibility}

\address Department of Mathematics, University of Kentucky,
Lexington, KY 40506.
\endaddress

\email shenz\@ms.uky.edu
\endemail

\thanks
The author is supported in part by the NSF (DMS-0500257).
\endthanks

\abstract 
Let $\OO$ be a bounded Lipschitz domain in $\br^n$.
We develop a new approach to
the invertibility on $L^p(\partial\OO)$ of the layer potentials
associated with elliptic equations and systems
in $\OO$.
As a consequence,
for $n\ge 4$ and $\frac{2(n-1)}{n+1}-\e<p<2$
where $\e>0$ depends on $\OO$, we obtain the solvability
of the $L^p$ Neumann type boundary value problems
for second order elliptic systems.
The analogous results for the biharmonic equation
are also established.
\endabstract



\endtopmatter

\document

\centerline{\bf 1. Introduction and Statement of Main Results}

Let $\OO$ be a bounded Lipschitz domain in $\br^n$.
The Dirichlet and Neumann problems 
for Laplace's equation in $\OO$
with boundary data in $L^p(\partial\OO)$
 had been well understood
more than twenty years ago.
Indeed it is known that the $L^p$ Dirichlet problem is uniquely solvable
for $2-\e<p\le \infty$, while the $L^p$ Neumann problem is
uniquely solvable for $1<p<2+\e$, where $\e>0$ depends on
$n$ and $\OO$.
Furthermore, the ranges of $p$'s are sharp; and
the solutions may be represented by
the classical layer potentials  \cite{D, JK, V1, DK1}.
Due to the lack of maximum principles and De Giogi -Nash
H\"older estimates,
the attempts to extend these results to second order elliptic
systems as well as to higher order elliptic equations had been successful
only in the case $n\ge 2$ for $p$ close to $2$ \cite{DKV1,
FKV, DKV2, F, K1, G, PV3, V2, V3}, 
and in the
lower dimensional case $n=2$ or $3$ for the sharp ranges of $p$'s
\cite{DK2, PV1, PV2, PV4}. 
Recently in \cite {S3, S4},
 we introduced a new approach to the $L^p$
Dirichlet problem
via $L^2$ estimates, reverse H\"older
inequalities and a real variable argument.
For second order elliptic systems as well as higher order elliptic
equations,
this led to the solvability of the $L^p$ Dirichlet problem
for $n\ge 4$ and $2<p<\frac{2(n-1)}{n-3}+\e$.
In the case of elliptic equations of order $2\ell$,
the upper bound of $p$ is known to be sharp for $4\le n\le 2\ell +1$ and
$\ell\ge 2$ \cite{PV3, PV4}.

The main purpose of this paper is to study the solvability
of the $L^p$ Neumann type boundary value problems
for elliptic systems and higher order equations.
We develop a new approach that can be used to establish the
$L^p$ invertibility of the trace operators $\pm (1/2)I+\Cal{K}^*$
of the double layer potentials for a limited range of $p$'s.
This limited-range approach is essential to the higher
order elliptic equations, as the $L^p$ invertibility
of $\pm (1/2)I +\Cal{K}^*$ fails in general for large $p$
in higher dimensions.
By duality, the invertibilty of $\pm (1/2)I+\Cal{K}^*$ on $L^{p}$
implies the invertibility of the Neumann trace operators
$\pm (1/2)I+\Cal{K}$ on $L^{p^\prime}$ of 
the single layer potentials. As a consequence, we are able to solve
the $L^p$ Neumann type problems for $p$
in the dual range $\frac{2(n-1)}{n+1}-\e_1
<p<2$.
We remark that in the lower dimensional case
$n=2$ or $3$, 
our approach recovers,
without the use of the Hardy spaces,
the $L^p$ solvability of the Neumann problem
for $1<p<2$ obtained in \cite{DK2} for elliptic systems.
The analogous results for the biharmonic equation, however,
are new even in the case $n=2$ or $3$.
It is also interesting to point out that
the approach we use here is in contrast with the method 
used in \cite{DK1}, where  the operators $\pm (1/2)I +\Cal{K}$
for the Neumann problem
are shown to be invertible first and 
the invertibility of $\pm (1/2)I +\Cal{K}^*$
for the Dirichlet problem is then established by duality.

This paper may be divided into three parts: elliptic systems,
the biharmonic equation, and Laplace's equation.
In the first part we consider  the
system of second order elliptic operators
$(\Cal{L}(\bu))^k
=-a_{ij}^{k\ell} D_iD_j u^\ell$ in  $\OO$, where
$D_i=\partial/\partial x_i$ and
$k,\ell=1,\dots, m$.
Let $N=(N_1, N_2,\dots, N_n)$ be the unit outward normal to
$\OO$ and
$$
\left(\frac{\partial \bu}{\partial\nu}\right)^k
=a_{ij}^{k\ell} \frac{\partial u^\ell}{\partial x_j} N_i
\tag 1.1
$$
denote the conormal derivatives of $\bu$ on $\partial\OO$.
We are interested in the $L^p$ Neumann type boundary value
 problem
$$
\left\{
\aligned
&\Cal{L}(\bu) =\bo\ \ \text{ in } \ \OO,\\
& \frac{\partial \bu}{\partial\nu} =\bbf
\in L^p(\partial\OO)\ \text{ on }\ \partial\OO,\\
& (\nabla \bu)^*\in L^p(\partial\OO),
\endaligned
\right.
\tag 1.2
$$
where $(\nabla\bu)^*$ denotes the nontangential maximal function
of $\nabla\bu$, and the boundary data $\bbf$ is taken in the sense
of nontangential convergence.
We will assume that $a_{ij}^{k\ell}$, $1\le i,j\le n$,
$1\le k,\ell \le m$ are real constants and satisfy the
symmetry condition $a_{ij}^{k\ell}
=a_{ji}^{\ell k}$ and the strong ellipticity condition
$$
\mu_0 \, |\xi|^2\le a_{ij}^{k\ell } \xi_i^k \xi_j^\ell
\le \frac{1}{\mu_0}\,
|\xi|^2,
\tag 1.3
$$
for some $\mu_0>0$ and any $\xi =(\xi_i^k)\in \br^{nm}$.
Let $\| \cdot\|_p$ denote the norm in $L^p(\partial\OO)$
with respect to the surface measure $d\sigma$ on $\partial\OO$.
The following is one of main results of the paper.

\proclaim{\bf Theorem 1.1}
Let $\OO$ be a bounded Lipschitz domain in $\br^n$, $n\ge 4$
with connected boundary. 
Then there exists $\e>0$ depending only on 
$n$, $m$, $\mu_0$ and $\OO$ such that, given
any $\bbf\in L^p(\partial\OO)$ with $\int_{\partial\OO}
\bbf \, d\sigma =\bo$ and
$$
\frac{2(n-1)}{n+1} -\e <p<2,
\tag 1.4
$$
the Neumann type problem (1.2) has a unique (up to constants)
solution $\bu$.
Furthermore, the solution $\bu$ satisfies the estimate
$\| (\nabla \bu)^*\|_p
\le C\, \| \bbf\|_p$ and may be represented by a single
layer potential with a density in $L^p(\partial\OO)$.
\endproclaim
 
Theorem 1.1 will be proved by the method of layer potentials.
Let $\Gamma (x)=(\Gamma^{k\ell}(x))_{m\times m}$ denote the matrix
of fundamental solutions for operator $\Cal{L}$ on $\br^n$. For
$\bg \in L^p(\partial\OO)$, let $\Cal{S}(\bg)$ and $\Cal{D}(\bg)$ denote
the single and double layer potentials respectively
with density $\bg$, defined by
$$
\align
& (\Cal{S}(\bg))^k(x)
=\int_{\partial\OO}
\Gamma^{k\ell}(y-x) \, g^\ell(y)\, d\sigma(y),
\tag 1.5\\
& (\Cal{D}(\bg))^k(x)
=\int_{\partial\OO}
\left\{ \frac{\partial}{\partial\nu(y)}
\Gamma_k (y-x)\right\}^\ell\, g^\ell(y)\, d\sigma(y),
\tag 1.6
\endalign
$$
where $\Gamma_k (x) =(\Gamma^{k1}(x), \dots, \Gamma^{km}(x))$
is the $k$th row of $\Gamma(x)$.
Let $\bu =\Cal{S}(\bg)$ and $\bv=\Cal{D}(\bg)$, then
$\Cal{L}(\bu)=\Cal{L}(\bv)=\bo$ in $\br^n\setminus
\partial\OO$. Moreover,
$$
\align
&\frac{\partial \bu_+}{\partial \nu}
=(\frac{1}{2} I +\Cal{K})\bg,\ \ \ 
\frac{\partial \bu_-}{\partial\nu}
=(-\frac12 I+\Cal{K})\bg,
\tag 1.7\\
&\bv_+ =(-\frac12 I +\Cal{K}^*) \bg, \ \ \ \bv_-=(\frac12 I +
\Cal{K}^*)\bg,
\tag 1.8
\endalign
$$
on $\partial\OO$, where $I$ denotes the identity operator, and
$\pm$ indicate the nontangential limits taken from
$\OO_+=\OO$ and $\OO_-=\br^n\setminus \overline{\OO}$
respectively. We remark that in (1.7)-(1.8),
 $\Cal{K}$ is a singular
integral operator on $\partial\OO$
and $\Cal{K}^*$ is the adjoint of $\Cal{K}$.
By \cite{CMM},
 $\Cal{K}$ and $\Cal{K}^*$ are bounded on $L^p(\partial\OO)$,
and $\| (\nabla \bu)^*\|_p +\| (\bv)^*\|_p
\le C\, \| \bg\|_p$
for any $1<p<\infty$.
In view of the trace formulas (1.7), the $L^p$ Neumann type problem
(1.2) is reduced to that of the invertibility of the operator
$(1/2) I +\Cal{K}$ on $L^p(\partial\OO)$ (modulo a finite dimensional
subspace). Similarly, because of (1.8),
one may solve the $L^p$ Dirichlet problem
$$
\left\{
\aligned
&\Cal{L}(\bu)=\bo \ \ \text{ in }\ \OO,\\
& \bu =\bbf \in L^p(\partial\OO)
\ \ \text{ on }\  \partial\OO,\\
&(\bu)^*\in L^p(\partial\OO),
\endaligned
\right.
\tag 1.9
$$
by showing that $-(1/2)I +\Cal{K}^*$ is invertible on $L^p(\partial\OO)$.
This is the so-called method of layer potentials for solving
boundary value problems.

For $n\ge 2$, the invertibility of $\pm (1/2)I +\Cal{K}$
on $L^p(\partial\OO)$ was indeed established in \cite{DKV2, FKV}
(also see \cite{K1, F, K2})
for $2-\e <p< 2+\e$, where $\e>0$ depends on the Lipschitz character of
$\OO$.
To do this, the main step is to show that for suitable solutions
of $\Cal{L}(\bu)=\bo$ in $\br^n\setminus \partial\OO$, one has
$$
\| \frac{\partial \bu_+}{\partial\nu}\|_2
\sim \| \nabla_t \bu_+\|_2
\ \ \ \text{ and }\ \ \ 
\| \frac{\partial \bu_-}{\partial\nu}\|_2+\| \bu\|_2
\sim \| \nabla_t \bu_-\|_2+\| \bu\|_2,
\tag 1.10
$$
where $\nabla_t \bu$ denotes the tangential derivatives of
$\bu$ on $\partial\OO$. As in the case of Laplace's equation \cite{V1},
the proof of (1.10) relies on the Rellich type
identities.

If we let $\bu=\Cal{S}(\bg)$ in (1.10), since
$\nabla_t \bu_+=\nabla_t \bu_-$ a.e. on $\partial\OO$, we obtain
$$
\| \frac{\partial \bu_+}{\partial \nu}\|_2 +\| \bu\|_2
\sim \| \frac{\partial\bu_-}{\partial\nu}\|_2 +\| \bu\|_2.
\tag 1.11
$$
It follows that
$$
\| \bg\|_2\le \| \frac{\partial \bu_+}{\partial\nu}\|_2
+ \| \frac{\partial \bu_-}{\partial\nu}\|_2
\le C\, \| (\pm \frac12 I+\Cal{K})\bg\|_2
+C\, \| \Cal{S}(\bg)\|_2.
\tag 1.12
$$
This is essentially enough to
deduce the invertibility of $\pm \frac12 I +\Cal{K}$ and hence
$\pm \frac12 I +\Cal{K}^*$ on $L^2(\partial\OO)$,
modulo some finite dimensional subspaces. By a perturbation argument
of A.P.~Calder\'on,
the invertibility can be extended to $L^p(\partial\OO)$ for
$p$ close to $2$. As a consequence,
the $L^p$ Dirichlet and Neumann type problems are solved
for $2-\e<p<2+\e$.

For Laplace's equation on Lipschitz domains,
 the invertibility of the corresponding
operators $\pm (1/2)I+\Cal{K}$ on $L^p(\partial\OO)$
was established for the sharp ranges of $p$'s in \cite{DK1}
(the case $p=2$ is in \cite{V1}).
The method used in \cite{DK1} relies on the classical
H\"older estimates for solutions of second order elliptic
equations of divergence form with bounded measurable coefficients.
Because of this,
the extension of the results in \cite{DK1} to elliptic
systems has only been successful
in the lower dimensional case $(n=2 \text{ or } 3)$ \cite{DK2}.
As we mentioned in the beginning of this section, we recently
introduced a new
approach to the $L^p$ Dirichlet problem for $p>2$ in \cite{S3, S4}.
Roughly speaking, this approach reduces the solvability of
the $L^p$ Dirichlet problem to a weak reverse H\"older
inequality on $I(P,r)$ with exponent $p$ for
$L^2$ solutions whose Dirichlet data vanish on $I(P,3r)$.
Here $I(P,r)=B(P,r)\cap\partial\OO$, where $P\in \partial\OO$
and $0<r<r_0$, is a surface ball
on $\partial\OO$. Combined with the $W^{1,2}$ regularity
estimate $\|(\nabla \bu)^*\|_2
\le C\, \|\nabla_t \bu\|_2$, this allows us
to establish the solvability of the
$L^p$ Dirichlet problem (1.9)
for $n\ge 4$ and
$$
2<p<\frac{2(n-1)}{n-3} +\e_1.
\tag 1.13
$$

In this paper we will show that if $\bv=\Cal{D}(\bg)$ is a double layer
potential, then 
$$
\| (\bv)^*\|_p \sim \| \bv_\pm\|_p,
\tag 1.14
$$
for any $p$ satisfying (1.13), where
the nontangential maximal function $(\bv)^*$
 is defined using nontangential approach
regions from both sides of $\partial\OO$.
Since $g=\bv_- - \bv_+$ by (1.8),
estimate (1.14) implies that $\pm (1/2)I +\Cal{K}^*$ are invertible on
$L^p(\partial\OO)$. By duality,
$\pm (1/2)I +\Cal{K}$ are invertible on $L^p(\partial\OO)$
for $p$ in the dual range (1.4).

By a refinement of the approach used in \cite{S3, S4},
we may reduce
the proof of (1.14) to the weak reverse
H\"older inequality
$$
\left\{\frac{1}{r^{n-1}}
\int_{I(P,r)}
|(\bv)^*|^p\, d\sigma\right\}^{1/p}
\le C\, \left\{
\frac{1}{r^{n-1}}
\int_{I(P,2r)}
|(\bv)^*|^2\, d\sigma\right\}^{1/2},
\tag 1.15
$$
where $\bv =\Cal{D}(\bg)$, and either $\bv_+ =\bo$ or 
$\bv_-=\bo$ on $I(P,3r)$.
The proof of (1.15) relies on applications of localized
$L^2$ estimates (or Rellich identities)
on the domains $B(P,r)\cap \OO_\pm$. It also depends on the fact
that 
$$
\frac{\partial \bv_+}{\partial\nu}
=\frac{\partial\bv_-}{\partial\nu}
 \ \ \ \ \text{ on }\ \partial\OO
\tag 1.16
$$
for any double layer potential $\bv$.
This crucial fact allows us to estimate 
the $L^2$ norm of $\nabla \bv_\pm$ on $I(P,r)$
by the $L^2$ norm of $\nabla_t \bv_\mp$ on $I(P,2r)$
respectively, plus some lower order terms.
See Lemma 2.4.
We mention that the upper bound of $p$ in (1.13) is dictated
by the use of Sobolev inequality on $I(P,r)$.
Whether this upper bound is necessary for the invertibility
of $\pm (1/2)I+\Cal{K}^*$ on $L^p(\partial\OO)$ for
second order elliptic systems remains open.

In this paper we also study the traction boundary value problem
for the system of elastostatics
$$
\left\{
\aligned
&\mu \Delta \bu +(\lambda +\mu)\nabla (\text{div}\, \bu)
=\bo\ \ \text{ in }\ \OO,\\
&\lambda (\text{div}\, \bu)N +\mu \big(\nabla \bu +(\nabla \bu)^T\big)N
=\bbf\in L^p(\partial\OO),\\
& (\nabla \bu)^*\in L^p(\partial\OO),
\endaligned
\right.
\tag 1.17
$$
where $\mu>0$, $ \lambda>-2\mu/n$ are Lam\'e constants, and $T$ indicates
the transpose of a matrix.
One may put (1.17) in the general form of (1.2) with
$$
a_{ij}^{k\ell}
=\mu \delta_{ij}\delta_{k\ell}
+\lambda \delta_{ik}\delta_{j\ell}
+\mu \delta_{i\ell}\delta_{jk}
\tag 1.18
$$
for $i,j,k,\ell =1,2,\dots, n$. It is easy to verify that
the coefficients satisfy the Legendre-Hadamard ellipticity
condition
$$
a_{ij}^{k\ell}\xi_i\xi_j \eta^k\eta^\ell\ge \mu \, |\xi|^2|\eta|^2
\ \ \text{ for any }\ 
\xi, \eta\in \br^n.
\tag 1.19
$$
However they do not satisfy the strong
elliptic condition (1.3). Thus Rellich type identities
alone are not strong enough to give estimate
(1.10). Nevertheless, this difficulty was overcome
in \cite {DKV2} by establishing a
 Korn type inequality on $\partial\OO$.
Consequently,
the $L^p$ traction problem (1.17)
was solved in \cite{DKV2} for $|p-2|<\e$. In the case $n=2$ or $3$,
the problem was solved in \cite{DK2} for the
optimal range $1<p<2+\e$.
Here we will show that
with a few modifications, the proof of Theorem 1.1 may be used to
solve the $L^p$ traction problem
for $p$ in the same range given in (1.4).
More specifically,
let $\Psi$ denote the space of vector valued functions 
$\bg =(g^1,\dots, g^n)$ on $\br^n$
satisfying
$
D_i g^j +D_j  g^i =0 \ \ \text{ for }\ 1\le i,j\le n$.
It is easy to show that $\bg\in \Psi$ if and only if
$\bg (x) =A x +\bold{b}$, where $\bold{b}\in \br^n$ and
$A$ is a real skew-symmetric matrix, $A^T=-A$. 
 Let
$$
L^p_\Psi
(\partial\OO)
=\big\{ \bbf\in L^p(\partial\OO):\ 
\int_{\partial\OO} \bbf \cdot \bg \, d\sigma
=0\  \ \text{ for all } \bg \in \Psi \big\}.
\tag 1.20
$$
 
\proclaim{\bf Theorem 1.2}
Let $\OO$ be a bounded Lipschitz domain in $\br^n$,
$n\ge 4$ with connected boundary.
Then there exists $\e>0$ depending only on $\lambda$,
$\mu$, $n$ and $\OO$ such that for
any $\bbf\in L^p_\Psi(\partial\OO)$ with $p$ satisfying (1.4),
the traction problem (1.17) has a solution $\bu$,
unique up to elements of $\Psi$.
Furthermore, the solution $\bu$ satisfies 
the estimate $\| (\nabla \bu)^*\|_p
\le C\, \| \bbf\|_p$ and may be represented by
a single layer potential with a density in $L^p(\partial\OO)$.
\endproclaim

The general program we outlined above for the second order
systems should apply to higher order elliptic
equations and systems, once the $L^2$ invertibility of
the layer potentials is established.
 In the second part of this paper, we study
the biharmonic Neumann problem
$$
\left\{
\aligned
 \Delta^2 u &=0 \ \ \text{ in }\ \OO,\\
\rho \Delta u +(1-\rho) \frac{\partial^2 u}{\partial N^2}
&=f\in L^p(\partial\OO)\ \ 
\text{ on }\ \partial\OO,\\
\frac{\partial}{\partial N} \Delta u
+\frac12 (1-\rho) \frac{\partial}
{\partial T_{ij}}
\left(\frac{\partial^2 u}{\partial N
\partial T_{ij}}\right)
&=\Lambda \in W^{-1,p}_0 (\partial\OO)
\ \ \text{ on }\ \partial\OO,
\endaligned
\right.
\tag 1.21
$$
where $\frac{\partial}{\partial T_{ij}}
=N_iD_j -N_j D_i$, and $W_0^{-1,p}(\partial\OO)$ denotes
the space of bounded linear functionals $\Lambda$ on 
$W^{1,p^\prime}(\partial\OO)$ such that
$\Lambda (1)=0$.
The $L^p$ Neumann problem (1.21) was recently formulated and studied
by G.~Verchota in \cite{V3}, where the solvability was established
for $p\in (2-\e, 2+\e)$ by the method of layer potentials.
 The following is the second main
result of the paper.

\proclaim{\bf Theorem 1.3}
Let $\OO$ be a bounded Lipschitz domain in $\br^n$, $n\ge 4$
with connected boundary. Let $(1/(1-n))<\rho<1$.
Then there exists $\e>0$ such that given any $f\in L^p(\partial\OO)$
and $\Lambda\in W^{-1,p}_0(\partial\OO)$ with
$\frac{2(n-1)}{n+1}-\e<p<2$, there exists a biharmonic
function $u$, unique up to linear functions, satisfying
(1.21) and $(\nabla\nabla u)^*\in L^p(\partial\OO)$.
Moreover, there exists a constant $C$ depending only on
$n$, $p$, $\rho$ and $\OO$ so that
$$
\| (\nabla\nabla u)^*\|_p\le C\, \big\{
\| \Lambda\|_{W^{-1,p}(\partial\OO)}
+\| f\|_p\big\},
\tag 1.22
$$
and the solution $u$ may be represented by a single layer potential.
If $n=2$ or $3$, above results hold for $1<p<2$.
\endproclaim

We refer the reader to Remark 7.3 for the ranges of $p$'s
for which the $L^p$ Dirichlet problem
for the biharmonic equation is uniquely solvable.
In particular the sharp ranges are known in the case $2\le n\le 7$.

In the last part of this paper we apply the method
used above for systems and the biharmonic equation
to the classical layer potentials for
Laplace's equation.
This allows us to recover the sharp $L^p$ results
in \cite{DK1}, without the use of the Hardy spaces.
In fact we are able to establish the following
stronger result.

\proclaim{\bf Theorem 1.4}
Let $\OO$ be a bounded Lipschitz domain  in  $\br^n$, $n\ge 3$
with connected boundary.
Then there exists $\delta>0$ depending only on $n$ and $\OO$, such
that
$$
\aligned
\frac12 I +\Cal{K}:&\ L^2_0\bigg(\partial \OO, 
\frac{d\sigma}{\omega}\bigg)
\to L^2_0\left(\partial\OO, \frac{d\sigma}{\omega}\right),\\
-\frac12 I +\Cal{K}^*: &\ L^2(\partial\OO,\omega d\sigma)
\to L^2(\OO,\omega d\sigma)
\endaligned
\tag 1.23
$$
are isomorphism for any $A_{1+\delta}$ weight $\omega$ on 
$\partial\OO$.
\endproclaim

We remark that the sharp $L^p$ invertibility of $(1/2)I +\Cal{K}$
and $-(1/2)I +\Cal{K^*}$ follows from Theorem 1.4
by an extrapolation theorem, due to Rubio de Francia \cite{R}.
Theorem 1.4 allows us to solve the Neumann problem for Laplace's equation
with boundary data in $L^2\left(\partial\OO,\frac{d\sigma}{\omega}\right)$.
This, combined with the weighted regularity estimate in \cite{S2},
shows that
$$
\| \frac{\partial u}{\partial N}\|_{L^2\left(\partial\OO,
\frac{d\sigma}{\omega}\right)}
\sim \| \nabla_t u\|_{L^2\left(\partial\OO,
\frac{d\sigma}{\omega}\right)},
\tag 1.24
$$
if $\Delta u=0$ in $\OO$ and $(\nabla u)^*\in L^2\left(\partial\OO,
\frac{d\sigma}{\omega}\right)$ with $\omega\in A_{1+\delta}(\partial\OO)$.

The paper is organized as follows.
Throughout Sections 2, 3 and 4,
 we will assume that the coefficients
$a_{ij}^{k\ell}$  of $\Cal{L}$
satisfy the symmetry condition 
$a_{ij}^{k\ell} =a_{ji}^{\ell k}$ and the
strong ellipticity condition (1.3).
In Section 2 we prove the reverse H\"older inequality
(1.15). See Theorem 2.6.
This is used in Section 3 to establish the invertibility of 
$\pm(1/2)I +\Cal{K}^*$ on $L^p$.
The proof of Theorem 1.1 is given in Section 4,
while the proof of
Theorem 1.2 can be found in Section 5.
Sections 6 and 7 deal with the biharmonic equation.
The corresponding reverse H\"older inequality
for biharmonic functions
is proved in section 6.
The proof of Theorem 1.3 is given in Section 7.
Finally the classical layer potentials are studied in
Section 8, where the proof of Theorem 1.4 can be found.
We point out that the usual conventions on 
repeated indices and on constants are used
throughout the paper.

\bigskip
 
\centerline{\bf 2. Reverse H\"older Inequalities}

Let $\OO$ be a bounded Lipschitz domain in $\br^n$.
Denote $\OO_+=\OO$ and $\OO_-=\br^n\setminus \overline{\OO}$.
For continuous function $u$ in $\OO_\pm$,
the nontangential maximal
function $(u)^*_\pm$ on $\partial\OO$ is defined by
$$
(\bu)^*_\pm (P) =\sup \big\{ |\bu(x)|: \
x\in \OO_\pm \ \ \text{ and }\ \
x\in \gamma(P)\big\},
\tag 2.1
$$
where $\gamma(P)=\{ x\in \br^n\setminus \partial\OO:\
|x-P|<2\, \text{dist}\, (x,\partial\OO)\}$. 

Assume
$0\in \partial\OO$ and
$$
\OO\cap B(0, r_0)
=\big\{ (x^\prime, x_n)\in \br^n:\ 
x_n>\psi(x^\prime)\big\} \cap B(0,r_0),
\tag 2.2
$$
where $\psi: \br^{n-1}\to \br$ is a Lipschitz function,
 and $\psi(0)=0$.
For $r>0$, we let
$$
 I_r 
=\big\{ (x^\prime, \psi(x^\prime))\in \br^{n-1}:\ 
|x_1|<r, \dots, |x_{n-1}|<r\big\},
\tag 2.3
$$
and
$$
\aligned
& D_r^+ =\big\{ (x^\prime, x_n):\ 
|x_1|<r,\dots, |x_{n-1}|<r,\
\psi(x^\prime)<x_n<\psi(x^\prime) +r\big\},\\
&D_r^- =\big\{ (x^\prime, x_n):\ 
|x_1|<r,\dots, |x_{n-1}|<r,\
\psi(x^\prime)-r<x_n<\psi(x^\prime)\big\}.
\endaligned
\tag 2.4
$$
Note that if $0<r<c\, r_0$, then $I_r\subset\partial\OO$
and $D_r^\pm\subset \OO_\pm$.

We begin with a boundary Cacciopoli's inequality.

\proclaim{\bf Lemma 2.1}
Suppose that $\Cal{L}(\bu)=\bo$ in $\OO_\pm$ and
$(\nabla\bu)^*_\pm \in L^2(I_{2r})$ for some $0<2r<c\, r_0$. Then
$$
\int_{D_r^\pm} |\nabla\bu|^2\, dx
\le \frac{C}{r^2}
\int_{D_{2r}^\pm} |\bu|^2\, dx
+C\, \int_{I_{2r}}\big|\frac{\partial \bu_\pm}{\partial \nu}\big|\,
|\bu_\pm|\, d\sigma.
\tag 2.5
$$
\endproclaim

\demo{Proof} The proof is rather standard. We first choose
a nonnegative function $\varphi \in C_0^\infty(\br^n)$ such that
$\varphi =1$ in $D^+_r$, $\varphi=0$ in $\OO\setminus D_{2r}^+$
and $|\nabla\varphi|\le C/r$. Let
$a(\xi,\eta)=a_{ij}^{k\ell}\xi_i^k \eta_j^\ell$ for
$\xi=(\xi_i^k),\ \eta=(\eta_j^\ell) \in \br^{mn}$.
It follows from integration by parts that
$$
\int_\OO a(\xi,\xi)\, \varphi^2\, dx
=-2\int_\OO a(\xi,\eta)\, \varphi\, dx
+\int_{\partial\OO} \frac{\partial\bu_+}{\partial\nu}\cdot
\bu_+\, \varphi^2\, d\sigma,
\tag 2.6
$$
where $\xi=(\xi_i^k)=(\frac{\partial u^k}{\partial x_i})$
and $\eta =(\eta_j^\ell)=(u^\ell\frac{\partial\varphi}{\partial x_j})$.
 Since $a(\xi,\xi)\ge 0$ for any $\xi\in \br^{mn}$, by Cauchy inequality,
we have
$$
|a(\xi,\eta)|\le a(\xi,\xi)^{1/2}\, a(\eta,\eta)^{1/2}
\le \frac14 a(\xi,\xi) + a(\eta,\eta).
\tag 2.7
$$
This, together with (2.6), gives
$$
\int_\OO a(\xi, \xi)\, \varphi^2\, dx
\le 4\int_\OO a(\eta, \eta)\, \varphi\, dx
+\int_{\partial\OO}
\frac{\partial\bu_+}{\partial\nu}\cdot
\bu_+\, \varphi^2\, d\sigma.
\tag 2.8
$$
Since $a(\xi, \xi)\ge \mu_0 |\nabla\bu|^2$, estimate (2.5) for the case
$D_r^+$ follows easily from (2.8). 
It is clear that the argument above also applies to the case $D_r^-$.
\enddemo

\proclaim{\bf Lemma 2.2}
Suppose that $\Cal{L}(\bu)=\bo$ in $\OO_\pm$ and
$(\nabla \bu)^*_\pm\in L^2(I_{2r})$ for some
$0<2r<c\, r_0$. Then
$$
\align
& \int_{I_r} |\nabla \bu_\pm|^2\, d\sigma
\le C\, \int_{I_{2r}}
\big|\frac{\partial \bu_\pm}{\partial\nu}\big|^2\, d\sigma
+\frac{C}{r}\int_{D^\pm_{2r}} |\nabla \bu|^2\, dx,\tag 2.9\\
& \int_{I_r} |\nabla \bu_\pm|^2\, d\sigma
\le C\, \int_{I_{2r}} |\nabla_t \bu_\pm|^2\, d\sigma
+\frac{C}{r}\int_{D^\pm_{2r}} |\nabla \bu|^2\, dx,\tag 2.10
\endalign
$$
where $\nabla_t \bu$ denotes the tangential derivatives of $\bu$ on
$\partial \OO$.
\endproclaim

\demo{Proof} To show (2.9), we observe that
 the $L^2$ Neumann
problem is solvable, uniquely up to constants,
on $D_{sr}^\pm$ for any $1<s<3/2$.
This yields
$$
\aligned
\int_{I_r}
|\nabla \bu_\pm |^2\, d\sigma
&\le \int_{\partial D^\pm_{sr}} |\nabla \bu|^2\, d\sigma
\le C\, \int_{\partial D_{sr}^\pm }
\big|\frac{\partial \bu}{\partial\nu}\big|^2\, d\sigma\\
&\le C\, \int_{I_{2r}}
\big|\frac{\partial \bu_\pm}{\partial\nu}\big|^2\, d\sigma
+C\, \int_{\OO_\pm \cap \partial D_{sr}^\pm}
|\nabla \bu|^2\, d\sigma.
\endaligned
\tag 2.11
$$
Estimate (2.9) now follows by
integrating both sides of (2.11) with respect to $s$ over 
interval $(1,3/2)$. Similarly, estimate (2.10) follows by applying
the regularity estimate 
$$
\int_{\partial D_{sr}^\pm}
|\nabla \bu|^2\, d\sigma
\le C\, \int_{\partial D_{sr}^\pm }
|\nabla_t \bu|^2\, d\sigma.
\tag 2.12
$$
for the Dirichlet problem on $D_{sr}^\pm$.
We remark that the regularity estimate
(2.12) and hence (2.10)
in fact hold for elliptic systems satisfying the
Legendre-Hadamard ellipticity condition (1.19) \cite{K1, F, G}.
This will be used in the proof of Theorem 1.2
\enddemo

In order to handle the solid integrals likes those in (2.9)-(2.10),
 we introduce a localized nontangential
maximal function,
$$
(\bu)^{*,r}_\pm(P) =\sup\big\{
|\bu(x)|:\ x\in \OO_\pm,\ |x-P|<c\, r \text{ and }\, 
|x-P|<2\, \text{dist}\, (x,\partial\OO)\big\}
\tag 2.13
$$
where $c>0$, depending on $\|\nabla\psi\|_\infty$ and $n$,
 is sufficiently small.

\proclaim{\bf Lemma 2.3}
Let $\bu$ be a continuous function on $D_{2r}^\pm$.
Then
$$
\left\{ \frac{1}{r^n}\int\Sb x\in D_r^\pm
\\ \delta(x)\le cr\endSb
|\bu|^p\, dx\right\}^{1/p}
\le C\left\{ \frac{1}{r^{n-1}}
\int_{I_{2r}}
|(\bu)^{*,r}_\pm |^q\, d\sigma\right\}^{1/q}
\tag 2.14
$$
where $\delta(x)=\text{dist}\, (x,\partial\OO)$ and
$1<q<p<nq/(n-1)$.
\endproclaim

\demo{Proof} We only consider the case $D_r^+$.
Note that
if $x=(x^\prime, x_n)\in D_r^+$ and $\delta(x)\le c\,r$, then
$|\bu(x)|\le (\bu)_+^{*,r}(y^\prime, \psi(y^\prime))$
for $|y^\prime -x^\prime|\le c\, \delta(x)$.
Hence, if $0<\alpha< n-1$,
$$
\aligned
|\bu(x)|\delta^\alpha (x)
&\le C\, \int_{|Q-P|<c\, \delta(x)}
\frac{(\bu)^{*,r}_+ (Q)}{|P-Q|^{n-1-\alpha}}\, d\sigma(Q)\\
&\le C\, \int_{|Q-P|<c\, r}
\frac{(\bu)^{*,r}_+(Q)}{|P-Q|^{n-1-\alpha}}\, d\sigma(Q),
\endaligned
\tag 2.15
$$
where $P=(x^\prime, \psi(x^\prime))$.
It follows that if $\alpha p<1$,
$$
\aligned
&\int\Sb x\in D_r^+\\ \delta(x)\le c\, r\endSb
|\bu(x)|^p\, dx\\
&\le C\, r^{1-\alpha p}
\int_{I_{r}} d\sigma(P)
\left\{
\int_{|Q-P|<c\, r}
\frac{(\bu)^{*,r}_+
(Q)}{|P-Q|^{n-1-\alpha}}\, d\sigma(Q)\right\}^p.
\endaligned
\tag 2.16
$$
This leads to the desired estimate (2.14) by the $L^q-L^p$ bounds of
the fractional integrals on $\partial\OO$ \cite{St1},
 where $1<q<p$ and 
$(1/q)-(1/p)=\alpha/(n-1)$. Finally we observe that
the condition $\alpha p<1$ is equivalent to $p<qn/(n-1)$. 
\enddemo

\proclaim{\bf Lemma 2.4}
Suppose that $\Cal{L}(\bu)=\bo$ in $\br^n\setminus\partial{\OO}$.
Assume that $\bu_+=\bo$ on $I_{32r}$
and $(\nabla \bu)_+^*
+(\nabla\bu)^*_-\in L^2(I_{32r})$ for some
$0<32r<c\, r_0$. Then
$$
\aligned
\int_{I_r} |\nabla \bu_-|^2\, d\sigma
&\le \frac{C}{r^2}
 \int_{I_{4r}} |\bu_-|^2\, d\sigma
+ \frac{C}{r^3}\int_{D_{32r}^+\cup D_{32r}^-}
|\bu|^2\, dx\\
&\ \ \ \ \ \ \ \ \ \ \ \ +C\, \int_{I_{4r}}
\big|\frac{\partial \bu_+}{\partial\nu}
-\frac{\partial \bu_-}{\partial\nu}\big|^2\, d\sigma.
\endaligned
\tag 2.17
$$
Similarly, if $\bu_-=\bo$ on $I_{32r}$, we have
$$
\aligned
\int_{I_r} |\nabla \bu_+|^2\, d\sigma
&\le 
\frac{C}{r^{2}}
 \int_{I_{4r}} |\bu_+|^2\, d\sigma
+ \frac{C}{r^3}\int_{D_{32r}^+\cup D_{32r}^-}
|\bu|^2\, dx\\
&\ \ \ \ \ \ \ \ \ \ \ \ \ +C\, \int_{I_{4r}}
\big|\frac{\partial \bu_+}{\partial\nu}
-\frac{\partial \bu_-}{\partial\nu}\big|^2\, d\sigma.
\endaligned
\tag 2.18
$$
\endproclaim

\demo{Proof}
Assume $\bu_+=\bo$ on $I_{32r}$. By using 
(2.9) and (2.5) as well as Cauchy inequality, we have
$$
\aligned
\int_{I_r} |\nabla \bu_-|^2\, d\sigma
& \le C\, \int_{I_{4r}}
\big|\frac{\partial \bu_-}{\partial\nu}\big|^2\, d\sigma
+\frac{C}{r^{2}}
 \int_{I_{8r}} |\bu_-|^2\, d\sigma\\
&\ \ \ \ \ \ \ \ \ \ \ \ \ \ \ \
+\frac{C}{r^3}
\int_{D_{4r}^-} |\bu|^2\, dx.
\endaligned
\tag 2.19
$$
Similarly, by (2.10) and (2.5), we obtain
$$
\int_{I_{4r}}
\big|\frac{\partial\bu_+}{\partial\nu}\big|^2\, d\sigma
\le  \frac{C}{r^3}
\int_{D_{32r}^+}|\bu|^2\, dx.
\tag 2.20
$$
where we have used the assumption $\bu_+=\bo$ and
hence $\nabla_t \bu_+=\bo$
on $I_{32r}$. Using 
$|\frac{\partial \bu_-}{\partial\nu}|
\le |\frac{\partial \bu_+}{\partial\nu}\big|
+|\frac{\partial \bu_+}{\partial\nu}
-\frac{\partial \bu_-}{\partial\nu}|$, it is not hard to see
that (2.17) follows from (2.19) and (2.20).
The proof of (2.18) is exactly the same.
\enddemo

Observe that estimates (2.17) and (2.18), together with the Sobolev
inequality
$$
\aligned
&\left\{ \frac{1}{|I_r|}\int_{I_r}
|\bu|^{{{p_n}}}\, d\sigma
\right\}^{1/{{p_n}}}\\
&\ \ \ \ \ \ 
\le C\,r \left\{ \frac{1}{|I_r|}
\int_{I_r} |\nabla_t \bu|^2\, d\sigma\right\}^{1/2}
+C\,
\left\{
\frac{1}{|I_r|} 
\int_{I_r} |\bu|^2\, d\sigma\right\}^{1/2},
\endaligned
\tag 2.21
$$
where $p_n=\frac{2(n-1)}{n-3}$ for $n\ge 4$, and
$p_3$ may be any exponent in $(2, \infty)$,
allows us to control the $L^{{p_n}}$ average of $\bu$ over $I_r$
by its $L^2$ average over $I_{4r}$, provided we can handle the last two
terms in the right sides of (2.17) and (2.18).
Since we will apply (2.17)-(2.18) to solutions given by the 
double layer potentials plus possible corrections,
the term involving $\frac{\partial \bu_+}{\partial\nu}
-\frac{\partial\bu_-}{\partial \nu}$ is negligible in view of
(1.16). In order to manage the remaining solid integrals,
it will be convenient to work with the nontangential maximal function
of $\bu$.

If $\bu$ is a function on $\br^n\setminus \partial\OO$, we let
$
(\bu)^* (P)=\max \{ (\bu)^*_+(P), \, (\bu)^*_-(P)\}
$
and
$$
 (\bu)^{*,r} (P)
=\sup\big\{ |\bu(x)|:\ \ x\in \gamma(P)\ \ 
\text{ and }\ \ |x-P|<c\, r\big\},
\tag 2.22
$$
for $P\in \partial\OO$,
where $c>0$ is sufficiently small.
By a simple geometric observation, we have
$$
\aligned
&\left\{ \frac{1}{|I_r|}
\int_{I_r}
|(\bu)^*|^p\, d\sigma\right\}^{1/p}\\
&\le 
\left\{ \frac{1}{|I_r|}
\int_{I_r}
| (\bu)^{*,r}|^p\, d\sigma\right\}^{1/p}
+
 \frac{C}{|I_{2r}|}
\int_{I_{2r}}
|(\bu)^{*}|\, d\sigma
\endaligned
\tag 2.23
$$
for any $p>1$. 

\proclaim{\bf Lemma 2.5} Let $\bar{p}>2$.
Suppose that the $L^{\bar{p}}$ Dirichlet problem for operator $\Cal{L}$ is
uniquely solvable for any bounded Lipschitz domain in $\br^n$.
Then for any $\frac{2(n-1)}{n}<p\le 2$,
$$
\aligned
&\left\{ \frac{1}{|I_r|}
\int_{I_r}
|(\bu)^*|^{\bar{p}}\, d\sigma\right\}^{1/{\bar{p}}}\\
&\le C\, \left\{ \frac{1}{|I_{4r}|}
\int_{I_{4r}}
\big(|\bu_+|
+|\bu_-|\big)^{\bar{p}}\, d\sigma\right\}^{1/{\bar{p}}}
+C\, \left\{ \frac{1}{|I_{4r}|}
\int_{I_{4r}}
|(\bu)^*|^p\, d\sigma\right\}^{1/p},
\endaligned
\tag 2.24
$$
where $\Cal{L}(\bu)=\bo$ in $\br^n\setminus \partial\OO$ and
$(\bu)^*\in L^{\bar{p}}(I_{4r})$.
\endproclaim

\demo{Proof} Since the $L^{\bar{p}}$ Dirichlet problem is solvable on
the Lipschitz domain $D^\pm_{sr}$, we have
$$
\int_{I_r}
| (\bu)^{*,r}|^{\bar{p}}\, d\sigma
\le C\, \int_{\partial D^+_{sr}} |\bu|^{\bar{p}}\, d\sigma
+\int_{\partial D^-_{sr}} |\bu|^{\bar{p}}\, d\sigma
\tag 2.25
$$
for $s\in (3/2,2)$.
It follows by an integration in $s$ over $(3/2,2)$ that
$$
\int_{I_r} |(\bu)^{*,r}|^{\bar{p}}\, d\sigma
\le C\, \int_{I_{2r}}
\big(|\bu_+| +|\bu_-|\big)^{\bar{p}}\, d\sigma
+\frac{C}{r}
\int_{D_{2r}^+\cup D_{2r}^-}
|\bu|^{\bar{p}}\, dx.
\tag 2.26
$$
This, together with estimates (2.23) and (2.14), yields that
$$
\aligned
&\left\{ \frac{1}{|I_r|}
\int_{I_r} |(\bu)^*|^{\bar{p}}\, d\sigma\right\}^{1/{\bar{p}}}\\
&\le C\,
\left\{ \frac{1}{|I_{3r}|}
\int_{I_{3r}}
\big(|\bu_+| +|\bu_-|\big)^{\bar{p}}\, d\sigma\right\}^{1/{\bar{p}}}
+C\, \left\{ \frac{1}{|I_{3r}|}
\int_{I_{3r}}
|(\bu)^*|^q\, d\sigma\right\}^{1/q}
\endaligned
\tag 2.27
$$
for any $q>(n-1)\bar{p}/n$. 
Since the $L^q$ Dirichlet problem for $\Cal{L}$
is also uniquely solvable for any $2\le q<p$,
it is not hard to see that
one may deduce estimate (2.24) for $p=2$ from (2.27) by using 
above argument repeatedly to decrease the exponent $q$ in (2.27)
to $2$. From here another application of the 
argument reduces the exponent from $2$ to any $q$
in $(2(n-1)/n, 2)$. 
\enddemo

Finally we are ready to state and prove the desired 
reverse H\"older inequality for elliptic systems.

\proclaim{\bf Theorem 2.6}
Suppose that $\Cal{L}(\bu)=\bo$ in $\br^n\setminus \partial\OO$
and $n\ge 4$.
Assume that either $\bu_+=\bo$ or $\bu_-=\bo$ on $I_{64r}$.
 Then,
if $(\nabla \bu)^*\in L^2(I_{64r})$ and $(\bu)^*\in L^{p_n}(I_{64r})$,
 we have
$$
\aligned
&\left\{\frac{1}{|I_r|}
\int_{I_r} |(\bu)^*|^{{p_n}}\, d\sigma\right\}^{1/{p_n}}\\
&\le 
C\, \left
\{\frac{1}{|I_{64r}|}
\int_{I_{64r}} |(\bu)^*|^2\, d\sigma\right\}^{1/2}
+
C\,r\, 
\left\{ \frac{1}{r^{n-1}}
\int_{I_{32r}}
\big|\frac{\partial \bu_+}{\partial\nu}
-\frac{\partial \bu_-}{\partial\nu}\big|^2\, d\sigma\right\}^{1/2},
\endaligned
\tag 2.28
$$
where $p_n=\frac{2(n-1)}{n-3}$.
If $n=3$, estimate (2.28) holds for any $p_3>2$.
\endproclaim

\demo{Proof} It is proved in \cite{S3} that
if $2<p<\frac{2(n-1)}{n-3}+\e$,
 the $L^p$ Dirichlet problem is uniquely solvable
for any bounded Lipschitz domain in $\br^n$. Thus
estimate (2.24) holds for $\bar{p}=p_n$.
This, combined with the Sobolev inequality (2.21), gives
$$
\aligned
&\left\{\frac{1}{|I_r|}
\int_{I_r} |(\bu)^*|^{{p_n}}\, d\sigma\right\}^{1/{p_n}}\\
&\le C\, r\, \left\{
\frac{1}{|I_{4r}|}
\int_{I_{4r}}
\big(|\nabla_t\bu_+|+|\nabla_t\bu_-|\big)^2\, d\sigma
\right\}^{1/2}
+\left\{ \frac{1}{|I_{4r}|}
\int_{I_{4r}}
|(\bu)^*|^2\, d\sigma\right\}^{1/2}.
\endaligned
\tag 2.29
$$
We now use (2.17)-(2.18) to estimate the term in (2.29)
with the tangential
derivatives.
Note that the solid integrals in (2.17)-(2.18)
are easily bounded by the maximal function $(\bu)^*$.
 Estimate (2.28) then follows.
\enddemo

\bigskip

\centerline{\bf 3. Invertibility of Double Layer Potentials in $L^p$}

Given $\bg \in L^p(\partial\OO)$ for some $1<p<\infty$. Let
$\bu=\Cal{D}(\bg)$ be the double layer potential defined  in (1.6).
Then $\bu_+ =(-(1/2)I +\Cal{K}^*)\bg$ and
 $\bu_-=((1/2)I +\Cal{K}^*)\bg$ on $
\partial \OO$. Moreover, we have $(\nabla \bu)^*
\in L^p(\partial\OO)$ and
$\frac{\partial \bu_+}{\partial\nu}
=\frac{\partial \bu_-}{\partial \nu}$ on $\partial\OO$, if
$\nabla_t \bg\in L^p(\partial\OO)$.

Since $\OO_-$ is connected, 
the kernel of operator $(1/2)I +\Cal{K}$ on $L^2(\partial
\OO)$ is of dimension $m$. Suppose $\{\bbf_\ell, \ell=1,\dots, m\}$
 spans
the kernel. Then $\int_{\partial\OO} \bbf_\ell \, d\sigma\neq
\bo$, and $\Cal{S}(\bbf_\ell)$ is a nonzero constant vector
in $\overline{\OO}$.
Let
$$
\Cal{X}^p(\partial\OO)=\big\{ \bbf \in L^p(\partial\OO):\ \ 
\int_{\partial\OO} \bbf \cdot \bbf_\ell \, d\sigma =0,\ 
\text{ for all } \ell=1,\dots,m\big\}
\tag 3.1
$$
for $p\ge 2$.
Since $\Cal{S}: L^p(\partial\OO)\to W^{1,p}(\partial\OO)$
is invertible for some $p>2$ \cite{G}, $\bbf_\ell\in L^p(\partial\OO)$
for some $p>2$. Thus the  space $\Cal{X}^p$ is also
well defined for $p>2-\e$.
 It was proved in \cite{DKV2} that
$$
\aligned
\frac12 I +\Cal{K}^*:&\ \ \Cal{X}^p(\partial\OO)\to \Cal{X}^p(\partial\OO),\\
 -\frac12 I +\Cal{K}^*:& \ \ L^p(\partial\OO)\to L^p(\partial\OO)
\endaligned
\tag 3.2
$$
are isomorphisms if $n\ge 3$ and $|p-2|<\e$. 
In the case $n=3$, the operators in (3.2) are isomorphisms for
$2-\e<p<\infty$ \cite{DK2}.
The goal of this section is to
establish the invertibility of $\pm (1/2)I +\Cal{K}^*$ for
$n\ge 4$ and $2<p<(2(n-1)/(n-3)) +\e$.

\proclaim{\bf Theorem 3.1}
There exists $\e>0$, depending on $n$, $m$, $\mu$ and
the Lipschitz character of $\OO$, such that
the operators $\pm (1/2)I +\Cal{K}^*$
in (3.2) are isomorphisms for
$n\ge 4$ and $2<p<\frac{2(n-1)}{n-3}+\e$.
\endproclaim

The proof of Theorem 3.1 is based on a real variable argument,
 inspired by a paper of Caffarelli and Peral 
\cite{CP} (see also \cite{W}). In \cite{S3, S4}, the argument
 was used to
solve the $L^p$ Dirichlet problem for elliptic systems and higher order
elliptic equations. 
This real variable argument may be considered 
as a dual and refined version of the celebrated Calder\'on-Zygmund
Lemma. We should mention that a similar argument
with a different motivation
was also used in \cite{ACDH} (see also \cite{A}).

The real variable argument may be formulated as follows.

\proclaim{\bf Theorem 3.2}
Let $Q_0$ be a cube in $\br^{n}$ and $F\in L^1(2Q_0)$. Let
$p>1$ and $f\in L^q(2Q_0)$ for some $1<q<p$.
 Suppose that for each dyadic subcube $Q$
of $Q_0$ with $|Q|\le \beta |Q_0|$, there
exist two integrable functions $F_Q$ and $R_Q$ on $2Q$ such that
$|F|\le |F_Q|+|R_Q|$ on $2Q$, and
$$
\align
\left\{\frac{1}{|2Q|}
\int_{2Q} |R_Q |^p\, dx\right\}^{1/p}
& \le 
C_1\left\{ \frac{1}{|\alpha Q|}
\int_{\alpha Q}| F|\, dx
+\sup\Sb Q^\prime\supset Q\endSb
\frac{1}{|Q^\prime|}\int_{Q^\prime}
|f|\, dx\right\},\tag 3.3\\
\frac{1}{|2Q|}
\int_{2Q}|F_Q|\, dx
&\le C_2
\sup\Sb  Q^\prime\supset Q\endSb
\frac{1}{|Q^\prime|}\int_{Q^\prime}
|f|\, dx,\tag 3.4
\endalign
$$
where $C_1, C_2>0$ and $0<\beta<1<\alpha$.
Then 
$$
\left\{ \frac{1}{|Q_0|}
\int_{Q_0}
|F|^q\, dx\right\}^{1/q}
\le \frac{C}{|2Q_0|}
\int_{2Q_0} |F|\, dx
+C\, \left\{ \frac{1}{|2Q_0|}
\int_{2Q_0}|f|^q\, dx\right\}^{1/q},
\tag 3.5
$$
where
$C>0$ is a constant depending only on 
$p$, $q$, $C_1$, $C_2$, $\alpha$, $\beta$ and $n$.
\endproclaim

We postpone the proof of Theorem 3.2 to the end of this section.

\remark{\bf Remark 3.3} Because of the local nature of Theorem 3.2,
it may be extended easily to each coordinate patch of $\partial\OO$.
Indeed, assume that $0\in \partial\OO$ and
$\OO\cap B(0,r_0)$ is given by (2.2).
Consider the map $\Phi: \partial D=\big\{
(x^\prime, \psi(x^\prime)):\ x^\prime\in \br^{n-1} \big\}
\to \br^{n-1}$, defined by $\Phi(x^\prime, \psi(x^\prime))
=x^\prime$.
We say $Q\subset \partial D$ is a surface cube of $\partial D$ if
$\Phi(Q)$ is a cube of $\br^{n-1}$. Moreover, a dilation
of $Q$ may be defined by $\alpha Q
=\Phi^{-1}(\alpha \Phi(Q))$. With these notations, 
one may state the extension of Theorem 3.2 to $\partial D$
in exactly the same manner as for the case of $\br^{n-1}$. Of course
in the case of $\partial D$, the constant $C$ in (3.5) also depends
on $\| \nabla \psi\|_\infty$.
\endremark

\demo{\bf Proof of Theorem 3.1}
We will give the proof for the invertibility of
 $(1/2)I +\Cal{K}^*$ on $\Cal{X}^p(\partial\OO)$. The case of
$-(1/2)I +\Cal{K}^*$ on $L^p(\partial\OO)$
 is similar and slightly easier.

Let $\bbf\in \Cal{X}^p(\partial\OO)\cap W^{1,2}(\partial\OO)$ for some
$p>2$. Since $(1/2)I +\Cal{K}^*$ is invertible
on $\Cal{X}^2(\partial\OO)$ and on $W^{1,2}(\partial\OO)/
\text{span}\{ \bbf_1, \dots, \bbf_m\}$, 
there exists $\bg\in \Cal{X}^2(\partial\OO)\cap W^{1,2}(\partial\OO)$
 such that
$((1/2)I+\Cal{K}^*)\bg =\bbf$ and $\| \bg\|_2\le C\, \| \bbf\|_2$.
Let $\bu=\Cal{D}(\bg)$ in $\br^n\setminus \partial\OO$. We will show that
there exists $\e>0$, depending only on 
$n$, $m$, $\mu_0$ and $\OO$, such that if $2<p<p_n+\e$,
$$
\aligned
&\left\{ \frac{1}{s^{n-1}}
\int_{B(P,s)\cap\partial\OO}
| (\bu)^*|^p\, d\sigma\right\}^{1/p}\\
&\le
C\,\left\{ \frac{1}{s^{n-1}}
\int_{B(P,Cs)\cap \partial\OO}
| (\bu)^*|^2\, d\sigma\right\}^{1/2}
+C\, \left\{\frac{1}{s^{n-1}}
\int_{B(P,Cs)\cap\partial\OO}
|\bbf|^p\, d\sigma\right\}^{1/p},
\endaligned
\tag 3.6
$$
for any $P\in \partial\OO$ and $s>0$ small.
Since $|\bg|=|\bu_+-\bu_-|\le 2\,(\bu)^*$, by covering
$\partial\OO$ with a finite number of small balls, estimate (3.6)
implies that
$$
\| \bg\|_p\le C\, \| \bg\|_2 +C\, \| \bbf\|_p
\le C\, \| \bbf\|_p.
\tag 3.7
$$
This shows that $(1/2)I +\Cal{K}^*: \, \Cal{X}^p
(\partial\OO) \to \Cal{X}^p(\partial\OO)$
is invertible, since $\Cal{X}^p(\partial\OO)\cap W^{1,2}(\partial\OO)$
is dense in $\Cal{X}^p(\partial\OO)$.

To prove (3.6), we use Theorems 3.2 and 2.6. By translation and rotation,
we may assume that $P=0$ and $B(0,r_0)\cap\OO$ is given
by (2.2). We consider the surface cube $Q_0=I_s$, defined in (2.3)
for $0<s<c\, r_0$. Let $Q$ be a small subcube of $Q_0$.
Choose $\varphi\in C_0^1(\br^n)$ such that
$\varphi=1$ on $200Q$, $\varphi=0$ in $\partial\OO
\setminus 300Q$ and
$|\nabla\varphi|\le C/r$, where $r$ is the diameter of $Q$. 
Since $L^2(\partial\OO)=\Cal{X}^2(\partial\OO) \oplus \br^m$, there exist
$\bg_Q\in \Cal{X}^2(\partial\OO)\cap W^{1,2}(\partial\OO)$ and $\bold{b}
\in \br^m$ such that
$$
\bbf \varphi
=(\frac12 I +\Cal{K}^*) \bg_Q +\bold{b}\ \ \ \text{ on }\ \partial\OO,
\tag 3.8
$$
and 
$
\| \bbf\varphi\|_2\sim \| \bg_Q\|_2
+|\bold{b}|
$.
Let
$
\bv =\Cal{D}(\bg_Q) +\bold{b}$
 in $ \br^n\setminus \partial\OO$ 
and $\bw=\bu-\bv$. 

We will apply Theorem 3.2 with $F=|(\bu)^*|^2$, $f=|\bbf|^2$ and 
$$
F_Q=2 |(\bv)^*|^2 \ \ \ \text{ and }\ \ \ R_Q=2|(\bw)^*|^2.
\tag 3.9
$$
Note that by the $L^2$ estimates,
$$\aligned
\frac{1}{|2Q|}
\int_{2Q}|F_Q|\, d\sigma
&\le \frac{C}{|Q|}
\int_{\partial\OO}
|(\bv)^*|^2\, d\sigma
\le \frac{C}{|Q|}
\left\{ \| \bg_Q\|_2^2 +|\bold{b}|^2\right\}\\
&\le \frac{C}{|200Q|}
\int_{200Q}
|\bbf|^2\, d\sigma.
\endaligned
\tag 3.10
$$
This gives condition (3.4). To verify (3.3), we observe that
$\bw_-=\bu_--\bv_-= \bbf (1-\varphi)$ on $\partial\OO$.
Hence $\bw_-=\bo$ on $200Q$. Also note that $(\nabla \bw)^*
\in L^2(\partial \OO)$ since $\bg$, $\bg_Q\in W^{1,2}(\partial\OO)$.
 It follows that
$(\bw)^*\in L^{p_n}(\partial\OO)$ (see e.g. \cite{S1}, p.1094).
Since $\bw=\Cal{D}(\bg)-\Cal{D}(\bg_Q)-\bold{b}$, we have
$ \frac{\partial \bw_+}{\partial\nu}
=\frac{\partial \bw_-}{\partial \nu}$ on $\partial\OO$.
Thus we may apply Theorem 2.6 to obtain
$$
\left\{ \frac{1}{|Q^\prime|}\int_{Q^\prime}
|(\bw)^*|^{p_n}\, d\sigma\right\}^{1/p_n}
\le C\,\left\{\frac{1}{|64Q^\prime|}\int_{64Q^\prime}
|(\bw)^*|^2\, d\sigma\right\}^{1/2},
\tag 3.11
$$
where $Q^\prime$ is any subcube of $Q$.
It is well known that the reverse H\"older inequalities like
(3.11) have the self-improving property (see e.g. \cite{Gi}).
This implies that there exists $\e>0$,
 depending only on $n$, $\|\nabla\psi\|_\infty$
and the constant $C$ in (3.11), such that
$$
\left\{ \frac{1}{|Q|}\int_{Q}
|(\bw)^*|^{\bar{p}}\, d\sigma\right\}^{1/\bar{p}}
\le C\,\left\{\frac{1}{|2Q|}\int_{2Q}
|(\bw)^*|^2\, d\sigma\right\}^{1/2}
\tag 3.12
$$
where $\bar{p}=p_n +\e$.
The right side of (3.12) may be estimated 
using $(\bw)^*\le (\bu)^* +(\bv)^*$
and then (3.10).
Thus condition (3.3) in Theorem 3.2 holds for
$p=p_n+\e$.
Consequently, 
estimate (3.6) holds for $2<p<p_n+\e$.
The proof is complete.
\enddemo
 
We now give the proof of Theorem 3.2.
The argument
 is essentially the same as that in the proof of
Lemma 2.18 in \cite{S3}. We shall need a localized
Hardy-Littlewood maximal function
$$
M_{Q} (g)(x)
=\sup\Sb Q^\prime\owns x\\
Q^\prime\subset Q \endSb
\frac{1}{|Q^\prime|}
\int_{Q^\prime} |g|\, dx
\tag 3.13
$$
for $x\in Q$, where $Q^\prime$ is a subcube of $Q$.

\demo{\bf Proof of Theorem 3.2}
For $\lambda>0$, let
$$
E(\lambda)=\big\{ x\in Q_0:\ 
M_{2Q_0}(F)(x)>\lambda\big\}.
\tag 3.14
$$
We claim that for any $1<q<p$,
it is possible to choose three constants
 $0<\delta<1$, $\gamma>0$ and $C_0>0$ depending
only on $n$, $C_1$, $C_2$, $\alpha$, $\beta$ in (3.3)-(3.4)
and $p,q$ such that
$$
|E(A\lambda)|\le \delta |E(\lambda)|
+|\big\{ x\in Q_0:\ M_{2Q_0}(f)(x)>\gamma\lambda\big\}|
\tag 3.15
$$
for all $\lambda>\lambda_0$, where $A=(2\delta)^{-1/q}$ and
$$
\lambda_0=\frac{C_0}{|2Q_0|}
\int_{2Q_0} |F|\, dx.
\tag 3.16
$$
Multiplying both sides of (3.15) by $\lambda^{q-1}$ and
then integrating the resulting inequality in 
$\lambda\in (\lambda_0, \Lambda)$, we obtain
$$
\int_{\lambda_0}^\Lambda \lambda^{q-1} |E(A\lambda)|\,
d\lambda
\le \delta \int_{\lambda_0}^\Lambda 
\lambda^{q-1} |E(\lambda)|\, d\lambda
+C_\gamma \, \int_{2Q_0} |f|^q\, dx,
\tag 3.17
$$
where we have used the fact that $M_{2Q_0}$ is bounded on $L^q$.
By a change of variable in the left side of (3.17), we
may deduce that
$$
A^{-q}(1-\delta A^{q})
\int_0^\Lambda \lambda^{q-1} |E(\lambda)|\, d\lambda
\le C\, |Q_0|\lambda_0^q
+C_\gamma \, \int_{2Q_0} |f|^q\, dx.
\tag 3.18
$$
Note that $\delta A^q
=1/2<1$.
Let $\Lambda\to \infty$ in (3.18). This gives
$$
\int_{Q_0} |F|^q\, dx
\le C\, |Q_0|\lambda_0^q +C\, \int_{2Q_0}|f|^q\, dx,
\tag 3.19
$$
which is (3.5) in view of (3.16).

To prove (3.15), we first note that $|E(\lambda)|\le C_n |Q_0|/C_0$
for any $\lambda>\lambda_0$.
This follows from the weak $(1,1)$ estimate for $M_{2Q_0}$.
Thus we may choose $C_0=2C_n/\delta$ so that $|E(\lambda)|<\delta\,|Q_0|$
for any $\lambda>\lambda_0$. 
We now fix $\lambda>\lambda_0$.
Since $E(\lambda)$ is open relative
to $Q_0$, we may write $E(\lambda)=\bigcup_k Q_k$, where
$Q_k$ are maximal dyadic subcubes of $Q_0$
 contained in $E(\lambda)$.
By choosing $\delta$ sufficiently small, we may certainly assume that
$|Q_k|<\beta |Q_0|$ and
$(\alpha+64)Q_k\subset 2Q_0$.

We will show that it is possible to choose $\delta>0$ and
$\gamma>0$ so that
$$
|E(A\lambda)\cap Q_k|\le \delta |Q_k|,
\tag 3.20
$$ 
whenever $\{ x\in Q_k:\ M_{2Q_0}(f)(x)\le \gamma \lambda\}
\neq \emptyset$. Clearly, estimate (3.15) follows from (3.20)
by summation.

Let $Q_k$ be such a maximal dyadic subcube.
Observe that
$$
M_{2Q_0}(F)(x)
\le \max \big\{ M_{2Q_k}(F)(x), C_n\lambda\big\},
\tag 3.21
$$
for any $x\in Q_k$. This is because $Q_k$ is maximal and so
$$
\frac{1}{|Q^\prime|}
\int_{Q^\prime} |F|\, dx
\le C_n\, \lambda
\tag 3.22
$$
for any $Q^\prime\cap Q_k\neq \emptyset$ and
$|Q^\prime|\ge c_n |Q_k|$.
We may assume that $A>C_n$. Then
$$
\aligned
 |E(A\lambda)\cap Q_k|&\le
|\big\{ x\in Q_k:\ M_{2Q_k} (F)>A\lambda\big\}|\\
&\le |\big\{ x\in Q_k:\ M_{2Q_k}(F_{Q_k})(x)>\frac{A\lambda}{2}\big\}|\\
&\ \ \ \ \ \ +
|\big\{ x\in Q_k:\ M_{2Q_k}(R_{Q_k})(x)>\frac{A\lambda}{2}\big\}|\\
&\le \frac{C_n}{A\lambda}
\int_{2Q_k}
|F_{Q_k}|\, dx
+\frac{C_{n,p}}{(A\lambda)^p}
\int_{2Q_k} |R_{Q_k}|^p\, dx,
\endaligned
\tag 3.23
$$
where we have used $|F|\le |F_{Q_k}| +|R_{Q_k}|$ on $2Q_k$ as well as weak
$(1,1)$, weak $(p,p)$ bounds of $M_{2Q_k}$.

By assumption (3.4), we have
$$
\aligned
\int_{2Q_k}|F_{Q_k}|\, dx&\le C_2\, |2Q_k|
\sup\Sb 2Q_0\supset Q^\prime
\supset Q_k\endSb
\frac{1}{|Q^\prime|}
\int_{Q^\prime} |f|\, dx\\
&\le C_2\, |2Q_k|\cdot \gamma \lambda,
\endaligned
\tag 3.24
$$
where the last inequality follows from the
fact $\{ x\in Q_k:\ M_{2Q_0}(f)\le \gamma\lambda\}
\neq \emptyset$. Similarly, we may use (3.3) and (3.22)
to obtain
$$
\aligned
\int_{2Q_k} |R_{Q_k}|^p\, dx
&\le C_1^p\cdot |2Q_k|
\left\{\frac{1}{|\alpha Q_k|}
\int_{\alpha Q_k} |F|\, dx
+\gamma \lambda\right\}^p\\
&
\le C_{n,\alpha}\, C_1^p\, |Q_k|\, \big\{ \lambda +\gamma \lambda\big\}^p.
\endaligned
\tag 3.25
$$
We now use (3.24) and (3.25) to estimate the right side
of (2.23). This yields 
$$
\aligned
|E(A\lambda)\cap Q_k|
&\le |Q_k|\left\{ \frac{C_n\, C_2\,  \gamma}{A}
+\frac{C_{n,\alpha, p}\, C_1^p }{A^p}\right\}\\
&
=\delta\, |Q_k|
\big\{ C_n \, C_2 \, \gamma \, \delta^{-\frac{1}{q}-1}
+C_{n,p,\alpha} \, C_1^p\,  \delta^{\frac{p}{q}-1}\big\}.
\endaligned
\tag 3.26
$$
Finally we observe that since $q<p$, it
is possible to choose
 $\delta>0$ so small that 
$$
C_{n,p,\alpha}\, C_1^p\,
\delta^{\frac{p}{q}-1}<(1/4).
$$
 After $\delta$ is chosen, we 
then choose $\gamma>0$ so small that $C_n\, C_2\, \gamma \, 
\delta^{-\frac{1}{q}
-1}<1/4$. This finishes the proof of (3.20) and thus the theorem.
\enddemo

The following weighted version of Theorem 3.2 will be used in Section 8.

\proclaim{\bf Theorem 3.4}
Under the same assumption as in Theorem 3.2, we have
$$
\left\{ \frac{1}{\omega(Q_0)}
\int_{Q_0}
|F|^q\, \omega dx\right\}^{1/q}
\le \frac{C}{|2Q_0|}
\int_{2Q_0} |F|\, dx
+C \left\{ \frac{1}{\omega(2Q_0)}
\int_{2Q_0} |f|^q\, \omega dx\right\}^{1/q},
\tag 3.27
$$
where $\omega$ is an $A_q$ weight on $2Q_0$ with the property that
for some $\eta>q/p$, 
$$
\frac{\omega(E)}{\omega(Q)}
\le C\, \left(\frac{|E|}{|Q|}\right)^\eta,
\tag 3.28
$$
 for any $E\subset Q\subset Q_0$.
\endproclaim

\demo{Proof}
Fix $1<q<p$. Since $\eta>q/p$,
we may choose $q_1\in (q,p)$ so that $\eta>q/q_1$.
Let $A=(2\delta)^{-1/q_1}$ in the proof of Theorem 3.2.
Note that if $|E(A\lambda)\cap Q_k|\le \delta |Q_k|$, then
$ \omega(E(A\lambda)\cap Q_k)\le C\, \delta^\eta \omega (Q_k)$.
This follows from (3.28). Thus
$$
\omega(E(A\lambda))
\le C\, \delta^\eta \omega(E(\lambda))
+\omega \left\{ x\in Q_0: \, M_{2Q_0} (f)>\gamma\lambda\right\},
\tag 3.29
$$
for any $\lambda\ge \lambda_0$. We now multiply both sides of (3.29)
by $\lambda^{q-1}$ and integrate the resulting inequality in
$\lambda$ from $\lambda_0 $ to $\Lambda$.
By a change of variable, we obtain
$$
\aligned
(A^{-q} -C\delta^\eta)
\int_0^\Lambda \lambda^{q-1} \omega(E(\lambda))\, d\lambda
&\le C\, \lambda_0^{q}\, \omega(Q_0)
+C_\delta \, \int_{Q_0} \big| M_{2Q_0}(f)|^q\, \omega\, dx\\
&\le C \, \lambda_0^q \, \omega(Q_0)
+C_\delta \, \int_{2Q_0} |f|^q\, \omega\, dx,
\endaligned
\tag 3.30
$$
where the second inequality follows from the well known
property of $M_{2Q_0}$ on $L^q(2Q_0, \omega\, dx)$
with $A_q$ weigh $\omega$ (see e.g. \cite{St2}).
Finally we note that since $\delta>q/q_1$,
we have $A^{-q} -C\delta^\eta
=(2\delta)^{q/q_1}-C\, \delta^\eta>0$ if $\delta>0$ is sufficiently
small. Estimate (3.27) follows from
(3.30) by letting $\Lambda\to\infty$.
\enddemo

\remark{\bf Remark 3.5}
If condition (3.3) holds for any $1<p<\infty$ (constant $C_1$ may
depend on $p$), then estimate (3.27) in Theorem 3.4 holds
for any $\omega\in A_q$.
This is because $w\in A_q$ implies condition (3.28) for some
$\eta=\eta(\omega)>0$.
\endremark

\bigskip

\centerline{\bf 4. The $L^p$ Boundary Value Problems for Elliptic Systems}
 
In this section we give the proof of Theorem 1.1 
stated in the Introduction. 
Let 
$$
L_0^p(\partial\OO)
=\big\{ \bbf \in L^p(\partial\OO):\ \ 
\int_{\partial\OO} \bbf\, d\sigma =\bo\big\}.
\tag 4.1
$$

\proclaim{\bf Theorem 4.1}
There exists $\e_1>0$, depending on $n$, $m$, $\mu_0$, and
the Lipschitz character of $\OO$, such that operators
$(1/2) I+\Cal{K}: L^p_0(\partial\OO)\to L^p_0(\partial\OO)$
and $-(1/2) I+\Cal{K}: L^p(\partial\OO)
\to L^p(\partial\OO)$ are invertible
for $\frac{2(n-1)}{n+1}-\e_1 <p<2$.
\endproclaim

\demo{Proof}
Let $p_0=\frac{2(n-1)}{n-3}+\e$, where $\e>0$ is
 given in Theorem 3.1. Note that $p_0^\prime
<\frac{2(n-1)}{n+1}$.
Since $-(1/2) I+\Cal{K}^*:\ L^p(\partial\OO)
\to L^p(\partial\OO)$ is invertible for $2<p<p_0$, by
duality,
we see that $-(1/2) I +\Cal{K}: L^p(\partial\OO)
\to L^p(\partial\OO)$ is invertible for
$p_0^\prime<p<2$.

Let $\bbf\in L^p_0(\partial\OO)$ for some
$p_0^\prime<p<2$. Given any $\bg\in L^{p^\prime}(\partial\OO)$,
since $L^{p^\prime}(\partial\OO)
=\Cal{X}^{p^\prime}(\partial\OO)\oplus\br^m$ and
$(1/2) I+\Cal{K}^*$ is invertible on $\Cal{X}^{p^\prime}(\partial\OO)$
by Theorem 3.1,
there exist $\bold{h}\in \Cal{X}^{p^\prime}(\partial\OO)$ and
$\bold{b}\in \br^m$ such that
$\bg =((1/2) I +\Cal{K}^*) \bold{h} +\bold{b}$
and $\| \bg\|_{p^\prime}\sim \| \bold{h}\|_{p^\prime}
+|\bold{b}|$.
Thus
$$
\aligned
\left|
\int_{\partial\OO} \bbf\cdot \bg\, d\sigma\right|
&=\left|\int_{\partial \OO}
\big(\frac12 I +\Cal{K}\big)\bbf\cdot \bold{h}\, d\sigma\right|\\
&\le \| \big(\frac12 I+\Cal{K}\big)\bbf\|_p\,  \| \bold{h}\|_{p^\prime}
\le C\, \|\big(\frac12 I +\Cal{K}\big) \bbf\|_p
\, \| \bg\|_{p^\prime}.
\endaligned
\tag  4.2
$$
It follows by duality that $\| \bbf\|_p
\le C\, \|((1/2) I +\Cal{K})\bbf\|_p$ for
any $\bbf \in L^p_0(\partial\OO)$. This shows
that $(1/2) I+\Cal{K}: \ L^p_0(\partial\OO)
\to L^p_0(\partial\OO)$ is one-to-one
and the range is closed. Note that the range is also
dense in $L^p_0(\partial\OO)$. This is because 
the operator is known to be invertible on $L^2_0(\partial\OO)$. 
Thus we have proved that
$(1/2) I+\Cal{K}$ is invertible on $L^p_0(\partial\OO)$
for any $p_0^\prime <p<2$.
\enddemo

\demo{\bf Proof of Theorem 1.1}
The existence follows directly from the invertibility of $(1/2) I+\Cal{K}$
on $L^p_0(\partial\OO)$ for $\frac{2(n-1)}{n+1}-\e_1
<p<2$.

In order to prove the uniqueness, 
we construct a matrix of the Neumann functions
$$
G_\nu^x(y) =\Gamma(x-y) -W^x(y),
\tag 4.3
$$
where for each $x\in \OO$, $W^x$ is a matrix solution
of the $L^2$ Neumann problem (1.2) with boundary data
$$
\frac{\partial}{\partial\nu(y)}
\big\{ \Gamma(x-y)\big\} +\frac{1}{|\partial\OO|} I_{m\times m}.
\tag 4.4
$$
In (4.4), $I_{m\times m}$ denotes the $m\times m$ identity
matrix.
By the $L^{2+\e}$ estimates for the Neumann problem,
we have $(\nabla W^x)^*\in L^p(\partial\OO)$ for some
$p>2$. Consequently, $(W^x)^*\in L^{p_1}(\partial\OO)$
for some $p_1>\frac{2(n-1)}{n-3}$ (see \cite{S1}, p.1094).

Suppose  now that $\Cal{L}(\bu)
=\bo$ in $\OO$, $(\nabla \bu)^*\in L^p(\partial\OO)$
and $\frac{\partial\bu}{\partial\nu} =\bo$ on $\partial\OO$.
Note that if $p>\max(p_0^\prime, p_1^\prime)$, then
$(\nabla \bu)^* (W^x)^* \in L^1(\partial\OO)$.
Similarly, one may show that $(\bu)^* (\nabla W^x)^*
\in L^1(\partial\OO)$. Thus one can use the integration by
parts, justified by the Lebesgue dominated convergence
theorem, to obtain the representation formula
$$
\aligned
\bu(x)&=
\int_{\partial\OO} G^x_\nu (y)\, \frac{\partial\bu}
{\partial\nu}\, d\sigma(y)
-\int_{\partial \OO} \frac{\partial G_\nu^x}{\partial \nu}
\bu(y)\, d\sigma(y)\\
&=\frac{1}{|\partial\OO|}
\int_{\partial\OO} \bu\, d\sigma.
\endaligned
\tag 4.5
$$
Hence $\bu$ is constant in $\OO$. 
The proof is finished.
\enddemo

\remark{\bf Remark 4.2}
Theorem 1.1 also holds in the exterior domain $\OO_-=\br^n\setminus
\overline{\OO}$ if one imposes additional condition
$|\bu(x)|=O(|x|^{n-2})$ as $|x|\to\infty$.
In this case the mean zero condition on $\bbf$ is not needed.
The proof is similar.
\endremark

\remark{\bf Remark 4.3}
Since $-(1/2) I+\Cal{K}^*$ is invertible on $L^p(\partial\OO)$
for $2<p<\frac{2(n-1)}{n-3}+\e$, the unique solution
of the $L^p$ Dirichlet problem (1.9), which was solved in 
\cite{S3}, may be represented by the double layer potential
$$
\bu(x)=\Cal{D}\big((-\frac12 I+\Cal{K}^*)^{-1}(\bbf)\big)(x).
\tag 4.6
$$
Since $L^p(\partial\OO)
=\Cal{X}^p(\partial\OO)\oplus \br^m$,
in the case of $\OO_-$, the solution may be represented
as $\bu =\Cal{D}(\bg) +\Cal{S}(\bold{h})$,
where $\bg\in \Cal{X}^p(\partial\OO)$, $\bold{h}
\in \text{Ker}((1/2) I +\Cal{K})$, and
$\| \bu\|_p\sim \| \bg\|_p +\| \bold{h}\|_p$.
\endremark

\remark{\bf Remark 4.4}
The Dirichlet problem with boundary data in $W^{1,p}(\partial\OO)$
for the elliptic systems satisfying the Legendre-Hadamard
condition (1.19) was solved in \cite{S3} for $n\ge 4$ and
$\frac{2(n-1)}{n+1}-\e<p<2$. This, combined with Theorem 1.1, gives
$
\|\frac{\partial \bu}{\partial \nu}\|_p
\sim \|\nabla_t \bu\|_p
$
for any solution of (1.2) with $p$ in the range (1.4).
\endremark
 
\bigskip

\centerline{\bf 5. The Traction Boundary Value Problem}

Throughout this section we assume that 
$$
\align
&\Cal{L}(\bu)=-\mu \Delta \bu -(\lambda +\mu)\nabla (\text{div}\,
\bu)\ \ \ \ \text{ in }\ \OO,
\tag 5.1\\
&\frac{\partial\bu}{\partial\nu}
=\lambda (\text{div}\,\bu) N +\mu \big(\nabla \bu
+(\nabla \bu)^T\big)N
\ \ \text{ on }\ \partial\OO.
\tag 5.2
\endalign
$$
If we write $(\Cal{L}(\bu))^k=-a_{ij}^{k\ell}D_iD_j u^\ell$, the conormal
derivatives (5.2) correspond to the choice of coefficients
given by (1.18). Note that $a_{ij}^{k\ell}$
do not satisfy the strong ellipticity condition (1.3).
However one has
$$
a_{ij}^{k\ell}\frac{\partial u^k}{\partial x_i}
\frac{\partial u^\ell}{\partial x_j}
=\lambda \, |\text{div}\, \bu|^2 
+\frac{\mu}{2}\, |\nabla\bu +(\nabla \bu)^T|^2
\sim  |\nabla\bu +(\nabla \bu)^T|^2.
\tag 5.3
$$ 
Using this observation, by establishing
 a Korn type inequality on the boundary,
Dahlberg, Kenig and Verchota were able to strength the Rellich
type inequalities.
This allows them to show that
$$
\aligned
\frac12 I +\Cal{K}:&\ \ \  L^p_\Psi (\partial\OO) \to L^p_\Psi(\partial\OO),\\
 -\frac12 I +\Cal{K}:& \ \ \ L^p(\partial\OO)\to L^p(\partial\OO),
\endaligned
\tag 5.4
$$
are invertible for $|p-2|<\e$ and $n\ge 2$ \cite{DKV2}, where
$L^p_\Psi(\partial\OO)$ is defined in (1.20).
In the case $n=2$ or $3$, it was proved in \cite{DK2}
that the operators in (5.4) are invertible for the optimal
range $1<p<2+\e$.
 The goal of this section is to
prove the following.

\proclaim{\bf Theorem 5.1}
There exists $\e>0$, depending on $n$, $\lambda$, $\mu$ and the
Lipschitz character of $\OO$, such that
the operators in (5.4) are invertible if $n\ge 4$ and
 $\frac{2(n-1)}{n+1}
-\e <p<2$.
\endproclaim

Let $\text{Ker}((1/2)I+\Cal{K})$ denote
 the kernel of operator $(1/2)I+\Cal{K}$ on 
$L^2(\partial\OO)$.
If $\bu=\Cal{S}(\bg)$ for some $\bg\in \text{Ker}((1/2)I+\Cal{K})$,
 then
$\frac{\partial \bu_+}{\partial\nu}=\bo$ on
$\partial\OO$. It follows from (5.3) and
integration by parts that $\nabla\bu +(\nabla \bu)^T=0$
in $\OO$. Thus $\Cal{S}(\bg)|_\OO\in \Psi$.
It is not hard to show that the map $\bg \to \Cal{S}(\bg)|_\OO$
from  Ker$((1/2)I+\Cal{K})$ to 
$\Psi$ is bijective.
Suppose $\{ \bg_k:\ k=1,2,\dots, n(n+1)/2\}$ spans 
$\text{Ker}((1/2)I+\Cal{K})$. Since $\Cal{S}:
L^p(\partial\OO)\to W^{1,p}(\partial\OO)$  is invertible for
$p$ close to $2$ \cite{G}, $\bg_k\in L^{q_0}(\partial\OO)$
for some $q_0>2$. Define
$$
\bold{T}^p(\partial\OO)
=\big\{ \bbf\in L^p(\partial\OO):\ 
\int_{\partial\OO} \bbf \cdot \bg_k\, d\sigma =0
\ \ \text{ for }k=1,2,\dots, n(n+1)/2\big\}
\tag 5.5
$$
for $p\ge q_0^\prime$.
 
\proclaim{\bf Theorem 5.2}
There exists $\e>0$ such that operators
$$
\aligned
\frac12 I +\Cal{K}^* & : \ \ \ \bold{T}^p(\partial\OO)\to 
\bold{T}^p(\partial\OO),\\
-\frac12 I +\Cal{K}^* & : \ \ \ L^p(\partial\OO) \to L^p(\partial\OO),
\endaligned
\tag 5.6
$$
are invertible for $n\ge 4$ and $2<p<\frac{2(n-1)}{n-3} +\e$.
\endproclaim

Theorem 5.1 follows from Theorem 5.2 by duality.
The case for $-(1/2)I +\Cal{K}$ is obvious. To see that
$(1/2)I +\Cal{K}$ is invertible on $L^{p}_\Psi(\partial\OO)$, we apply the
same duality argument as in the proof of Theorem 4.1. 
To do this, we only need to show that
$L^{p^\prime}(\partial\OO)=\bold{T}^{p^\prime}(\partial\OO)
\oplus \Psi$. By a dimensional consideration, it suffices to prove
that $\bold{T}^{p^\prime}(\partial\OO)\cap \Psi=\{ 0\}$.
To this end, let $\bg \in \bold{T}^{p^\prime}(\partial\OO)
\cap \Psi$. Then $\bg =\Cal{S}(\bold{h})$ on $\partial\OO$ for some
$\bold{h}\in \text{Ker}((1/2)I+\Cal{K})$.
Let $\bu=\Cal{S}(\bold{h})$ in $\br^n$. Since
$\bold{h}=\frac{\partial \bu_+}{\partial\nu}
-\frac{\partial \bu_-}{\partial\nu}=-\frac{\partial \bu_-}{\partial\nu}$,
we obtain
$$
\int_{\OO_-}
a_{ij}^{k\ell}\frac{\partial u^k}{\partial x_i}
\frac{\partial u^\ell}{\partial x_j}\, dx
=-\int_{\partial\OO}
\frac{\partial\bu_-}{\partial\nu}\cdot \bu\, d\sigma
=\int_{\partial\OO}
\bold{h}\cdot \bg\, d\sigma
=0,
\tag 5.7
$$
 where the last equality follows from the fact that
$\bg$ is in the range of $(1/2)I +\Cal{K}^*$ on $L^2(\partial\OO)$.
One may deduce from (5.7) that $\bu|_{\OO_-}\in \Psi$.
This implies that $\frac{\partial \bu_-}{\partial\nu} =\bo$ and
thus $\bold{h}=\bo$.

Since the proof of Theorem 5.2 uses the same line of argument
as in the proof of Theorem 3.1, we will only point out the
necessary modification needed here.

First, because of (5.3), estimate (2.5) is replaced by
$$
\int_{D_r^\pm} |\nabla\bu +(\nabla\bu)^T|^2\, dx
\le \frac{C}{r^2}
\int_{D_{2r}^\pm} |\bu|^2\, dx
+C\, \int_{I_{2r}}\big|\frac{\partial \bu_\pm}{\partial \nu}\big|\,
|\bu_\pm|\, d\sigma.
\tag 5.8
$$
The proof is exactly the same.

Next, estimate (2.9) needs to be modified, as we used
$$
\|\nabla \bu\|_{L^2(\partial D^\pm_{sr})}\le C\, \| \frac{\partial \bu}
{\partial\nu}\|_{L^2(\partial D^\pm_{sr})}
\tag 5.9
$$
for any $L^2$ solutions. In the case of (5.1), we know
that estimate (5.9) is 
true for one of such solutions, $\bv$, given by a single
layer potential with density
$((1/2)I+\Cal{K})^{-1}(\frac{\partial \bu}{\partial\nu})$.
 If $\bu$ is another solution with the same
traction boundary data on $\partial D_{sr}^\pm$, then $\bold{w}
=\bu-\bv =Ax +\bold{b}\in \Psi$. It follows that
$$
\aligned
\int_{\partial D_{sr}^\pm} |\nabla \bu|^2\, d\sigma
& \le C\, \int_{\partial D_{sr}^\pm} |\nabla \bold{v}|^2\, d
\sigma
+C\, r^{n-1}\ |A|^2\\
&
\le C\, \int_{\partial D_{sr}^\pm} \big|
\frac{\partial\bu}{\partial\nu}\big|^2\, d\sigma
+
C\, r^{n-1}\ |A|^2.
\endaligned
\tag 5.10
$$
Since $\bold{w}$ is a linear function and thus harmonic, we have
$$
\int_{D_{sr}^\pm} |\nabla \bold{w}|^2\, dx \le \int_{\partial D_{sr}^\pm}
|\bold{w}|\, |\nabla \bold{w}|\, d\sigma.
\tag 5.11
$$
It follows that
$$
\aligned
|A|&\le \frac{C}{r^n}\int_{\partial D_{sr}^\pm}
|\bold{w}|\, d\sigma
\le \frac{C}{r^n}\int_{\partial D_{sr}^\pm}
\big( |\bu| +|\bv|\big)\, d\sigma\\
&\le \frac{C}{r}
\left\{ \frac{1}{r^{n-1}}
\int_{\partial D_{sr}^\pm}
|\bu|^2\, d\sigma\right\}^{1/2}
+C\, \left\{ \frac{1}{r^{n-1}}
\int_{\partial D_{sr}^\pm}
\big|\frac{\partial \bu}{\partial\nu}\big|^2\, d
\sigma\right\}^{1/2}.
\endaligned
\tag 5.12
$$
This, together with (5.10), gives
$$
\int_{\partial D_{sr}^\pm}
|\nabla \bu|^2\, d\sigma
\le C\, \int_{\partial D_{sr}^\pm}
\big|\frac{\partial \bu}{\partial\nu}\big|^2\, d\sigma
+\frac{C}{r^2}\,
\int_{\partial D_{sr}^\pm}
|\bu|^2\, d\sigma.
\tag 5.13
$$
By integrating both sides of (5.13) in $s\in (1,3/2)$, we obtain
$$
\aligned
\int_{I_r}|\nabla\bu_\pm|^2\, d\sigma
&\le C\, \int_{I_{2r}}
\big|\frac{\partial\bu_\pm }{\partial\nu}\big|^2\, d\sigma
+\frac{C}{r}\int_{D_{2r}^\pm}
|\nabla\bu +(\nabla\bu)^T|^2\, dx\\
&\ \ \ \ \ \ \ \ \ \ \ \ \ \ \ \ \ 
+\frac{C}{r^{3}}
\int_{D^\pm_{2r}} |\bu|^2\, dx.
\endaligned
\tag 5.14
$$
This replaces estimate (2.9). The extra term in (5.14) is harmless.

Finally in the proof of Lemma 2.4, we used estimate (2.5) to
estimate the solid integral of $|\nabla\bu|^2$ on $D_{sr}^\pm$.
In the case of (5.1), we consider $\bv=\bu -Ax$, where
$$
A=\frac{1}{2|D_{sr}^\pm|}\int_{D_{sr}^\pm}
\big( \nabla\bu -(\nabla\bu)^T\big)\, dx.
\tag 5.15
$$
Then by Korn's inequality (see \cite{DKV2}, Lemma 1.18), we have
$$
\int_{D_{sr}^\pm}
|\nabla \bv|^2\, dx
\le C\, \int_{D_{sr}^\pm} |\nabla \bv +(\nabla \bv)^T|^2\, dx.
\tag 5.16
$$
Note that integration by parts gives
$$
|A|\le \frac{C}{r^n}\int_{\partial D_{sr}^\pm}
|\bu|\, d\sigma.
\tag 5.17
$$
It follows that
$$
\aligned
\int_{D_{sr}^\pm} |\nabla\bu|^2\,dx
&\le C\int_{D_{sr}^\pm}|\nabla\bu +(\nabla \bu)^T|^2\, dx
+C\, r^n\, |A|^2\\
&
\le C\int_{D_{sr}^\pm}|\nabla\bu +(\nabla \bu)^T|^2\, dx
+\frac{C}{r}\int_{\partial D_{sr}^\pm}
|\bu|^2\, d\sigma.
\endaligned
\tag 5.18
$$
We now integrate both sides of (5.18) in $s\in (1,3/2)$. This yields 
$$
\aligned
\int_{D_r^\pm} |\nabla\bu|^2\, dx
&\le C\, \int_{D^\pm_{2r}}
|\nabla \bu +(\nabla \bu)^T|^2\, dx
+\frac{C}{r}\int_{I_{2r}}
|\bu_\pm |^2\, d\sigma\\
&\ \ \ \ \ \ \ \ \ \ \ \ \
+\frac{C}{r^2}
\int_{D^\pm_{2r}}
|\bu|^2\, dx.
\endaligned
\tag 5.19
$$
Estimate (5.19), combined with (5.8), allows us to bound the solid
integral of $|\nabla \bu|^2$ in the same manner as in the 
strong elliptic case.
Because of this,  Lemma 2.4 and therefore
Theorem 2.6 hold for the
system of elastostatics.
 Consequently, Theorem 5.2 is proved using the same line of argument
as in the proof of Theorem 3.1.
We should point out that since $a_{ij}^{k\ell}$
satisfy the Legendre-Hadamard ellipticity condition,
the $L^p$ Dirichlet problem is solved for $2<p<\frac{2(n-1)}{n-3}
+\e$ and $n\ge 4$ in \cite{S3}.
This is used in the proof of Theorem 5.2.
We omit the details.

We end this section with 

\demo{\bf The Proof of Theorem 1.2}
The existence follows from the invertibility of $(1/2)I +\Cal{K}$ on $L^p_\Psi
(\partial\OO)$ for $p$ in the range given in (1.4). As in the case
of Theorem 1.1, to prove the uniqueness, 
one constructs a matrix Neumann function
$G_\nu^x (y)=\Gamma(x-y) -W^x(y)$, where $W^x$ is a matrix
whose $i$th row is an
$L^2$ solution of (1.17) with the traction boundary data
$$
\frac{\partial }{\partial \nu(y)}
\left\{ \Gamma_i (y-x)\right\} -\sum_{k=1}^{\frac{n(n+1)}{2}} C_{i,k}^x
\big\{ A_k y+\bold{b}_k\big\}.
\tag 5.20
$$
Here $\{ A_k y +\bold{b}_k,\, k=1,2,\dots, \frac{n(n+1)}{2}\}$
is an orthonormal basis of $\Psi$ with respect to the $L^2(\partial\OO)$
norm, and
$$
C_{i,k}^x  =\int_{\partial\OO} 
\frac{\partial }{\partial \nu(y)}
\left\{ \Gamma_i (y-x)\right\} \cdot (A_k y +\bold{b}_k)\,  d\sigma(y)
=-(A_k x +b_k)^i
\tag 5.21
$$
so that the functions in (5.20) belong to $L^2_\Psi(\partial\OO)$.
The same argument as in the proof of Theorem 1.1 shows
that if $\Cal{L}(\bu) =\bo$ in $\OO$, $(\nabla \bu)^*\in L^p(\partial
\OO)$ for some $p>\frac{2(n-1)}{n+1}-\e$, and $
\frac{\partial\bu}{\partial\nu}=\bo$ on $\partial\OO$, then
$$
\aligned
\bu (x)
&=-\int_{\partial\OO}\frac{\partial G^x_\nu}{\partial\nu}\, \bu\,
d\sigma\\
&
=(A_k x+\bold{b}_k)\int_{\partial \OO}
\{ A_k y +\bold{b}_k\} \cdot \bu(y)\, d\sigma(y).
\endaligned
\tag 5.22
$$
Thus $\bu\in \Psi$. This finishes the proof.
\enddemo

\bigskip

\centerline{\bf 6. Reverse H\"older Inequalities for
Biharmonic Functions}

For simplicity, we will assume that 
$\frac{1}{1-n}<\rho<1$. 
Some modifications are needed in the case $\rho=\frac{1}{1-n}$.
Following \cite{V3}, we let
$$
\aligned
 M_\rho (u) &=\rho \Delta u +(1-\rho) \frac{\partial^2 u}{
\partial N^2}=\rho\Delta u +(1-\rho) N_i N_j D_i D_j u,\\
 K_\rho (u)& =\frac{\partial\Delta u }{\partial N}
 +\frac12 (1-\rho) \frac{\partial}
{\partial T_{ij}}
\left(\frac{\partial^2 u}{\partial N\partial T_{ij}}\right)\\
&=\frac{\partial\Delta u }{\partial N}
+\frac12 (1-\rho) (N_iD_j -N_j D_i)
\big(N_k(N_i D_j-N_j D_i)D_k u\big),
\endaligned
\tag 6.1
$$
where $\frac{\partial}{\partial T_{ij}}
=N_i D_j -N_j D_i$.
Observe that $N_iN_j \frac{\partial u}{\partial T_{ij}}=0$.

Assume $0\in \partial\OO$ and $\OO\cap B(0,r_0)$ is given by
(2.2). Let $W^{1,2}(I_r)$ denote the space of functions $f$ on 
$I_r$ such that $|\nabla_t f|\in L^2(I_r)$, where $I_r$ is defined
in (2.3). We will use the scale-invariant norm
$$
\| f\|_{W^{1,2}(I_r)}
=\left\{ \int_{I_r} |\nabla_t f|^2\, d\sigma
+\frac{1}{r^2}
\int_{I_r} |f|^2\, d\sigma\right\}^{1/2}
\tag 6.2
$$
for $W^{1,2}(I_r)$, whose dual space is denoted by $W^{-1,2}(I_r)$.

The following is a boundary Cacciopoli inequality.

\proclaim{\bf Lemma 6.1} Suppose $\Delta^2 u=0$ in $\OO_\pm$ and
$(\nabla\nabla u)^*_\pm \in L^2(I_{3r})$. Then
$$
\aligned
\int_{D_r^\pm} |\nabla \nabla u|^2\, dx
&\le C\, \| u \varphi\|_{W^{1,2}(I_{2r})}
\, \| \varphi K_\rho(u
)\|_{W^{-1,2}
(I_{2r})}\\
&\ \ \ \ \ \ \ +C\, \| \frac{\partial (u \varphi^2)}{\partial N}\|_2
\, \| M_\rho(u )\|_{L^2(I_{2r})}\\
&\ \ \ \ \ \ \ +\frac{C}{r^2}\int_{D_{2r}^\pm}
|\nabla  u|^2\, dx
+\frac{C}{r^2}\int_{I_{2r}} |u |\, |\nabla u|\, d
\sigma,
\endaligned
\tag 6.3
$$
where $\varphi$ is a function in $C_0^\infty(B(0,(3/2)r))$ such that
$\varphi=1$ in $B(0,r)$, $0\le \varphi\le 1$ and $|\nabla\varphi|
\le C/r$.
\endproclaim

\demo{Proof}
Let $v=u\varphi^2$. It follows from the integration by parts and 
$\Delta^2 u=0$ in $\OO_\pm$ that
$$
\aligned
&\int_{\partial\OO}
\left\{ v\, K_\rho (u) -\frac{\partial v}{\partial N}
M_\rho (u)\right\}\, d\sigma\\
&=\mp\int_{\OO_\pm}
\big\{ (1-\rho) D_iD_j v
\cdot D_iD_j u
+\rho \Delta v\cdot \Delta u\big\}\, dx.
\endaligned
\tag 6.4
$$
Note that
$$
D_iD_j v\cdot D_iD_j u
=\varphi^2 |\nabla \nabla u|^2 +4\varphi D_iu D_j \varphi \cdot D_iD_j u
+uD_iD_j\varphi^2 \cdot D_iD_j u
\tag 6.5
$$
and $\Delta v \cdot \Delta u
=\varphi^2 |\Delta u|^2 +4\varphi D_iu D_i\varphi \cdot \Delta u
+u \Delta \varphi^2\cdot \Delta u$.
The second term in the right side of (6.5) can be absorbed
by the first term using the Cauchy inequality with an $\e$.
To handle the last term in the right side of (6.5), one uses
the integration by parts again. This produces the last integral in (6.3).
Finally, to finish the proof, we observe that
$$
(1-\rho) |\nabla\nabla u|^2 +\rho |\Delta u|^2\ge c_\rho
|\nabla \nabla u|^2,
\tag 6.6
$$
if $\frac{1}{1-n}<\rho<1$ (see \cite{V3}). 
\enddemo

\remark{\bf Remark 6.2} If, in addition, in Lemma 6.1
 we assume that
$u_\pm =|\nabla u_\pm|=0$ on $I_{2r}$, then
$$
\int_{D_r ^\pm }
|\nabla \nabla u|^2\, dx
\le \frac{C}{r^2}\, 
\int_{D_{2r}^\pm } |\nabla u|^2\, dx.
\tag 6.7
$$
This is the usual boundary 
Cacciopoli's inequality for the biharmonic equation.
\endremark

\remark{\bf Remark 6.3}
It follows from (6.3) and the Cauchy inequality with an $\e$ that
$$
\aligned
\int_{D_r^\pm}
|\nabla\nabla u |^2\, dx
&\le \e r\,  \| \varphi K_\rho(u)\|^2_{W^{-1,2}(\partial\OO)}
+\e r\, \| M_\rho (u)\|^2_{L^2(I_{2r})}\\
&\ \ \ \ \
+\frac{C_\e }{r}\int_{I_{2r}}
|\nabla u|^2\, d\sigma
+\frac{C}{r^2}\, \int_{D_{2r}^\pm} |\nabla u|^2\, dx.
\endaligned
\tag 6.8
$$
We remark that the intergals in (6.3)
which involve $|u|^2$ on $I_{2r}$
may be handled by replacing $|u|^2$ with $|u-c|^2$ and
using the Poincar\'e inequality.
\endremark

Our next lemma relies on the following Rellich type identity
discovered by G.~Verchota (\cite{V3}, pp.232-233) 
for the biharmonic equation,
$$
\aligned
&\frac12 \int_{\partial\OO}
<N,\alpha>\big\{ (1-\rho)|\nabla\nabla u|^2
+\rho |\Delta u|^2\big\}\, d\sigma\\
&\ \ \ =\int_{\partial\OO}
\frac{\partial}{\partial N} (\alpha\cdot \nabla u)
M_\rho(u)\, d\sigma
-\int_{\partial\OO}
(\alpha \cdot \nabla u)\, K_\rho(u)\, d\sigma\\
&\ \ \ \ \ \ \ \ \ \ \ \ \ \ \ \ \ \ \ 
\pm (1-\rho)\,\int_{\OO_\pm}
E_{ij}(\alpha, u)\, L_{ij}(u)\, dx,
\endaligned
\tag 6.9
$$
where $L_{ij}=D_iD_j +\theta \delta_{ij}\Delta$ and
$$
E_{ij}( \alpha, u)
=\frac12\,\text{div}(\alpha) \, L_{ij}(u)
-L_{ij}(\alpha)\cdot \nabla u
-2D_i \alpha\cdot \nabla D_j u
-2\theta \delta_{ij} D_k\alpha \cdot \nabla D_k u.
$$
In (6.9), $\alpha\in C_0^\infty(\br^n, \br^n)$ is a vector field and
$u$ is a suitable biharmonic function in $\OO_\pm$.
Also $\theta$ is related to $\rho$ by $\rho=(n\theta +n\theta^2)/(1+2\theta
+n\theta^2)$. With identity (6.9), Verchota was
able to extend the method of layer potentials from second order
equations and systems to the fourth order biharmonic equation.
This identity will also play a crucial role in our
study of the $L^p$ biharmonic Neumann problem.

\proclaim{\bf Lemma 6.4}
Under the same assumption as in Lemma 6.1, we have
$$
\aligned
\int_{I_r} |\nabla \nabla u|^2\, d\sigma
&\le C\, \| \varphi K_\rho(u)\|_{W^{-1,2}(I_{2r})}^2
+C\, \| M_\rho (u)\|^2_{L^2(I_{2r})}\\
&+ \frac{C}{r^2}\int_{I_{2r}}
|\nabla u|^2\, d\sigma
+\frac{C}{r}\int_{D^\pm_{2r}} |\nabla \nabla u|^2\, dx
+\frac{C}{r^3}
\int_{D^\pm_{2r}} |\nabla u|^2\, dx,
\endaligned
\tag 6.10
$$
where $\varphi\in C^\infty(B(0,(3/2)r))$ 
is the same function as in Lemma 6.1.
\endproclaim

\demo{Proof} Let $\alpha=-\bold{e}_n \varphi^2$ where
$\bold{e}_n=(0,\dots, 0,1)$. We apply the Rellich identity
(6.9) on the Lipschitz domain $D_{sr}^\pm$, where $s\in (3/2,2)$.
Since $<N,-\bold{e}_n>\ge c>0$ on $I_{2r}$,
this gives
$$
\aligned
&c\int_{I_{sr}}
|\varphi \nabla\nabla u|^2\, d\sigma\\
&\le C\, \int_{\OO_\pm \cap \partial D_{sr}^\pm}
|\nabla\nabla u|^2\, d\sigma
+C\, \| \varphi\nabla u\|_{W^{1,2}(I_{2r})}
\| \varphi K_\rho(u)\|_{W^{-1,2}(I_{2r})}\\
&\ \ \ \ \ \ \ \ \
+C\, \|\nabla (\alpha\cdot\nabla u)\|_{L^2(I_{2r})}
\| M_\rho(u)\|_{L^2(I_{2r})}\\
&\ \ \ \ \ \ \ \ \ \ 
+\frac{C}{r}\int_{D^\pm_{2r}}|\nabla \nabla u|^2\, dx
+\frac{C}{r^3}\,
\int_{D^\pm_{2r}} |\nabla u|^2\, dx.
\endaligned
\tag 6.11
$$
Using the Cauchy inequality with an $\e$, it is not hard to see that
 the higher order terms in $\| \varphi \nabla u\|_{W^{1,2}(I_{2r})}$ and
$\|\nabla (\alpha\cdot \nabla u)\|_{L^2(I_{2r})}$
 may be absorbed by
 the left side of (6.11).
Finally a familiar integration in $s$ over $(3/2,2)$
enables us to handle the first term in the right side of (6.10),
as in Section 2.
\enddemo

\remark{\bf Remark 6.5} Suppose $\Delta^2 u=0$ in $\OO_\pm$ and
$(\nabla \nabla u)^*_\pm \in L^2(I_{3r})$. If $u_\pm=|\nabla u_\pm|=0$
on $I_{2r}$, then
$$
\int_{I_r} |\nabla \nabla u|^2\, d\sigma
\le \frac{C}{r^3}
\int_{D_{2r}^\pm} |\nabla u|^2\, dx.
\tag 6.12
$$
This follows from the regularity estimate \cite{V2}
$$
\int_{\partial D_{sr}^\pm }
|\nabla \nabla u|^2\, d\sigma
\le C\, \int_{\partial D_{sr}^\pm}|\nabla_t \nabla u|^2
\, d\sigma,
\tag 6.13
$$
together with estimate (6.7), by an integration in $s\in (3/2,2)$.
\endremark

Recall that $(\nabla \nabla u)^* =\max \big\{
(\nabla \nabla u)^*_+, 
(\nabla\nabla u)^*_-\big\} $ for functions $u$ defined in $\br^n\setminus
\partial\OO$.

\proclaim{\bf Lemma 6.6}
Suppose $\Delta^2 u =0$ in $\br^n\setminus \partial\OO$ and
$(\nabla \nabla u)^*
\in L^2(I_{32r})$. Assume that either $u_+=|\nabla u_+|=0$ or
$u_-=|\nabla u_-|=0$ on $I_{32r}$. Then
$$
\aligned
\int_{I_r} |\nabla\nabla u_\pm|^2\, d\sigma
&\le \frac{C}{r^2}
\int_{I_{8r}}\big\{ |\nabla u_+|^2 +|\nabla u_-|^2\big\}\, d
\sigma
+\frac{C}{r^3}\, 
\int_{D_{16r}^+\cup D_{16r}^-}
|\nabla u|^2\, dx\\
&\ \ \ \ \ \ \ \ \ \ \ \ +C\, \| \varphi_1 \big[ K_\rho(u_+)-K_\rho (u_-)\big]
\|_{W^{-1,2}(I_{4r})}^2\\
&\ \ \ \ \ \ \ \ \ \ \ \ + 
C\, \| \varphi_2 \big[ K_\rho(u_+)-K_\rho (u_-)\big]
\|_{W^{-1,2}(I_{4r})}^2\\
& \ \ \ \ \ \ \ \ \ \ \ \ +
C\, \| M_\rho(u_+)-M_\rho(u_-)\|^2_{L^2(I_{4r})},
\endaligned
\tag 6.14
$$
where $\varphi_1$, $\varphi_2$ are two functions in $C^\infty_0(B(0,4r))$
with the properties that $0\le \varphi_i\le 1$ and
$|\nabla\varphi_i|\le C/r$ for $i=1,2$.
\endproclaim

\demo{Proof} Assume that $u_+=|\nabla u_+|=0$ on $I_{32r}$.
By (6.10) and (6.8), we obtain
$$
\aligned
\int_{I_r} |\nabla\nabla u_-|^2\, d\sigma
&\le C\, \| \varphi_1 K_\rho (u_-)\|^2_{W^{-1,2}(I_{4r})}
+C\, \| \varphi_2 K_\rho (u_-)\|^2_{W^{-1,2}(I_{4r})}\\
&\ \ \ \ \ \
+C\, \| M_\rho(u_-)\|^2_{L^2(I_{4r})}\\
&\ \ \ \ \ \
+\frac{C}{r^2}\int_{I_{4r}}|\nabla u_-|^2\, d\sigma
+\frac{C}{r^3}\,
\int_{D_{4r}^-}|\nabla u|^2\, dx
\endaligned
\tag 6.15
$$
where $\varphi_1\in C^\infty_0(B(0, (3/2)r))$ and
$\varphi_2\in C^\infty_0(B(0, 3r))$. In view of (6.14) and (6.15),
we need to estimate
$\|\varphi_i K_\rho(u_+)\|^2_{W^{-1,2}(I_{4r})}$, $i=1,2$ and
$\| M_\rho (u_+)\|_{L^2(I_{4r})}^2$. Clearly, by Remark 6.5,
$$
\| M_\rho(u_+)\|_{L^2(I_{4r})}^2
\le C\, \int_{I_{4r}}
|\nabla \nabla u_+|^2\, d\sigma
\le\frac{C}{r^3}
\int_{D_{8r}^+}
|\nabla u|^2\, dx.
\tag 6.16
$$

Finally, since supp$\varphi_i\subset B(0,3r)$,
the term $\|\varphi_i K_\rho(u_+)\|^2_{W^{-1,2}(I_{4r})}$ is bounded by
$$
\aligned
C\, \| \varphi_i &\frac{\partial}{\partial N}
\big(\Delta u_+\big)\|_{W^{-1,2}(I_{4r})}
+C\, \| \nabla\nabla u_+\|^2_{L^2(I_{4r})}\\
&
\le C\, \| \varphi_i \frac{\partial}{\partial N}
\big(\Delta u_+\big)\|^2_{W^{-1,2}(\partial D_{sr})}
+C\, \| \nabla\nabla u_+\|^2_{L^2(I_{4r})}\\
&\le C\, \| \Delta u_+\|^2_{L^2(\partial D_{sr})}
+C\, \| \nabla\nabla u_+\|^2_{L^2(I_{4r})}\\
& \le C\, \int_{I_{5r}}
|\nabla \nabla u_+|^2\, d\sigma
+C\, \int_{\OO\cap \partial D_{sr}^+}
|\nabla\nabla u|^2\, d\sigma,
\endaligned
\tag 6.17
$$
for any $s\in (4,5)$, where we have used the $L^2$ regularity estimate 
in $D_{sr}^+$ for
Laplace's equation in the second inequality. With (6.12) and (6.7)
at our disposal,
the desired estimate for
$\|\varphi_i K_\rho(u_+)\|^2_{W^{-1,2}(I_{4r})}$
now follows from (6.17) by an integration in $s\in (4,5)$.
The case $u_-=|\nabla u_-|=0$ on $I_{32 r}$ is exactly the same.
This completes the proof.
\enddemo

As in Section 2, estimate (6.14) leads to a reverse H\"older
inequality.

\proclaim{\bf Theorem 6.7} Under the same assumption as in Lemma
6.6, we have
$$
\aligned
\left\{ \frac{1}{|I_r|}\int_{I_r}
|(\nabla u)^*|^{p_n}\, d\sigma\right\}^{1/p_n}
& \le C\, 
\left\{ \frac{1}{|I_{32r}|}\int_{I_{32r}}
|(\nabla u)^*|^2\, d\sigma\right\}^{1/2}\\
&\ \ \ \ \ \ \ +C\, \| \varphi_1 \big[ K_\rho(u_+)-K_\rho (u_-)\big]
\|_{W^{-1,2}(I_{4r})}^2\\
&\ \ \ \ \ \ \ + 
C\, \| \varphi_2 \big[ K_\rho(u_+)-K_\rho (u_-)\big]
\|_{W^{-1,2}(I_{4r})}^2\\
& \ \ \ \ \ \ \ +
C\, \| M_\rho(u_+)-M_\rho(u_-)\|^2_{L^2(I_{4r})},
\endaligned
\tag 6.18
$$
where $p_n=\frac{2(n-1)}{n-3}$ for $n\ge 4$. If $n=2$ or $3$, estimate (6.18)
holds for any $2<p_n<\infty$.
\endproclaim

\demo{Proof}
The proof is similar to that of Theorem 2.6 with $\nabla u$ in the
place of $\bu$. We leave the details to the reader.
However we should remark that
the proof uses the solvability of the $L^{p_n}$
Dirichlet problem for the biharmonic equation on any bounded
Lipschitz domains. But this has been established in \cite{PV1} for $n=2$ or
$3$,
and in \cite{S3} for $n\ge 4$.
\enddemo

\bigskip

\centerline{\bf 7. The $L^p$ Biharmonic Neumann Problem}

This section is devoted to the proof of Theorem 1.3.
We begin with the definition of the
biharmonic layer potentials introduced by Verchota in \cite{V3}.
Fix $x\in \br^n$, let $B^x=B^x(y)$ denote the fundamental solution
for operator $\Delta^2$ with pole at $x$, given by
$$
B^x(y)=
\left\{
\alignedat 3
&\frac{1}{2(n-2)(n-4)\omega_n}\cdot \frac{1}{|x-y|^{n-4}},
&\quad\quad & n=3 \ \text{ or }\  n\ge 5,\\
&-\frac{1}{4\omega_4}\, \log |x-y|,
&\quad\quad & n=4,\\
&-\frac{1}{8\pi} |x-y|^2 \big( 1-\log |x-y|\big), &\quad\quad & n=2.
\endalignedat
\right.
\tag 7.1
$$
Given $(F,g)\in W^{1,p}(\partial\OO)\times L^p(\partial\OO)$
for $1<p<\infty$, the
double layer potential for the biharmonic equation is defined
by
$$
w(x)=\Cal{D}_\rho (F,g)(x)
=\int_{\partial\OO}
\big\{ K_\rho(B^x)(y) F(y) + M_\rho(B^x)(y) g(y)\big\}\,
d\sigma(y),
\tag 7.2
$$
for $x\in \br^n\setminus \partial\OO$. Clearly $\Delta^2 w=0$
in $\br^n\setminus \partial\OO$. By computing
$K_\rho(B^x)$ and $M_\rho(B^x)$ in (7.2), one may show
that
$$
w(x)=\int_{\partial\OO}
\left\{ \frac{\partial \Gamma^x}{\partial N}\, F
+\Gamma^x\, g
+(1-\rho)\frac{\partial }{\partial T_{jk}}
D_k B^x \cdot \left( 
N_i\frac{\partial F}{\partial T_{ij}}-N_j g\right)\right\}\,
d\sigma,
\tag 7.3
$$
where $\Gamma^x=\Delta B^x$ is the fundamental solution
for $\Delta$ with pole at $x$.
Also
$$
\aligned
 D_\ell w(x)
=& -\int_{\partial\OO}
\left\{
D_i \Gamma^x \cdot \frac{\partial F}{\partial T_{\ell i}}
+D_\ell\Gamma^x \cdot g\right\}\, d\sigma\\
& \ \ \  \ \ \ -(1-\rho)\, \int_{\partial\OO}
\left\{\frac{\partial }{\partial T_{jk}}
D_k D_\ell B^x \cdot \left( 
N_i\frac{\partial F}{\partial T_{ij}}-N_j g\right)\right\}\,
d\sigma.
\endaligned
\tag 7.4
$$
It follows by \cite{CMM} that
$$
\|(\nabla w)^*\|_p \le C\, \big\{
\| \nabla_t F\|_p
+\| g\|_p\big\}.
\tag 7.5
$$
To compute the nontangential limits of $w$ and $\nabla w$, one uses
$$
\aligned
\lim\Sb x\to P\in \partial\OO\\
x\in \OO_\pm\cap \gamma(P)\endSb &
\int_{\partial\OO}
D_iD_jD_k B^x\cdot f\, d\sigma\\
&=\pm \frac12 N_i N_j N_k f(P)
+\text{ p.v.}
\int_{\partial\OO}
D_iD_jD_k B^P\cdot f\, d\sigma.
\endaligned
\tag 7.6
$$
This, together with (7.3)-(7.4), gives
$$
\big(w_\pm ,-\frac{\partial w_\pm}{\partial N}\big)
=(\pm \frac12 +\Cal{K}_\rho^*) (F,g),
\tag 7.7
$$
where $\Cal{K}^*_\rho$ is a bounded operator on $W^{1,p}(\partial\OO)
\times L^p(\partial\OO)$.

For $(\Lambda, f)\in W^{-1,p}(\partial\OO)\times L^p(\partial\OO)$
with $1<p<\infty$,
the single layer potential is defined by
$$
v(x)=\Cal{S}(\Lambda, f)(x)
=\Lambda(B^x(\cdot)) -\int_{\partial\OO} \frac{\partial B^x}
{\partial N}\, f\, d\sigma.
\tag 7.8
$$
Clearly $\Delta^2 v=0$ in $\br^n\setminus \partial\OO$.
By writing
$
\Lambda =\frac{\partial h_{ij}}{\partial T_{ij}} +h_0
$ 
with $h_{ij}$, $h_0\in L^p(\partial\OO)$ so that
$$
\Lambda(B^x)=\int_{\partial\OO}
\left\{ -\frac{\partial B^x}{\partial T_{ij}}\, h_{ij}
+B^x\, h_0\right\}\, d\sigma,
\tag 7.9
$$
one sees that
$$
\| (\nabla \nabla v)^*\|_p\le C\, \big\{
\|\Lambda\|_{W^{-1,p}(\partial\OO)}
+\| f\|_p\big\}
\tag 7.10
$$
for $1<p<\infty$ by \cite{CMM}. Also
$$
\big(K_\rho(v)_\pm , M_\rho(v)_\pm\big)
=(\mp\frac12 I +\Cal{K}_\rho)(\Lambda, f)
\tag 7.11
$$
where operator $\Cal{K}_\rho$, whose adjoint is $\Cal{K}^*_\rho$ in (7.4),
is bounded on $W^{-1,p}(\partial\OO)\times L^p(\partial\OO)$.
We point out that the trace of $K_\rho(v)_\pm $ in (7.11)
 is taken in the sense of distribution, i.e.,
$$
K_\rho(v)_\pm (\phi)
=\lim_{k\to\infty}
\int_{\partial \OO_k^\pm }
K_\rho (v)\, \phi\, d\sigma,
\tag 7.12
$$
for $\phi\in C_0^1(\br^n)$, where $\OO_k^\pm$ is a sequence of 
smooth domains which approximate $\OO_\pm$ from inside,
respectively \cite{V1}.
Because of (7.6), to prove (7.11), we only need to take
care of the term $\frac{\partial}{\partial N}\Delta v$.
To do this, we note that
$$
\aligned
\Delta v& =-\int_{\partial\OO}
\left\{ \frac{\partial \Gamma^x}{\partial T_{ij}}\,
h_{ij} +\frac{\partial \Gamma^x}{\partial N}\, f\right\}\,
d\sigma
+\int_{\partial\OO}\Gamma^x\, h_0\, d\sigma\\
&=
D_j \int_{\partial\OO}
\Gamma^x \big\{ N_i h_{ij} -N_i h_{ji}
+N_j f\big\}\, d\sigma
+\int_{\partial\OO} \Gamma^x\, h_0\, d\sigma.
\endaligned 
\tag 7.13
$$
This allows us to express $\frac{\partial }{\partial N}
\Delta v$ on $\partial \OO_k$ in terms of tangential derivatives
plus a higher order term,
$$
\frac{\partial \Delta v}{\partial N}
=\frac{\partial}{\partial T_{\ell j}} 
D_\ell 
 \int_{\partial\OO}
\Gamma^x \big\{ N_i h_{ij} -N_i h_{ji}
+N_j f\big\}\, d\sigma
+\frac{\partial}{\partial N}
\int_{\partial\OO} \Gamma^x\, h_0\, d\sigma.
\tag 7.14
$$
We remark that the computation of the trace operators 
in \cite{V3} used the harmonic extension of functions in
$W^{1,p^\prime}(\partial\OO)$ to $\OO$.
On general Lipschitz domains, this would require $p>2-\e$.

Let $\bold{X}^p(\partial\OO)$ denote the subspace of
$W^{-1,p}(\partial\OO)\times L^p(\partial\OO)$ whose elements
$(\Lambda, f)$ satisfy
$$
\Lambda (1)=0\ \ 
\text{ and }\ 
\Lambda(x_j)=\int_{\partial\OO} f\, N_j\, d\sigma
\text { for } \ \ j=1,\dots, n.
\tag 7.15
$$
One of the main results in \cite{V3} is that
$$
\aligned
\frac12 I +\Cal{K}_\rho&: \ \ 
W^{-1,p}(\partial\OO)\times L^p(\partial\OO)
\to W^{-1,p}(\partial\OO)\times L^p(\partial\OO),\\
-\frac12 I+\Cal{K}_\rho&:\ \ 
\bold{X}^p(\partial\OO)
\to \bold{X}^p(\partial\OO)
\endaligned
\tag 7.16
$$
are isomorphism for $p\in (2-\e, 2+\e)$. Let
$\{ (\Lambda_j^*, f_j^*):\ j=0,1,\dots,n\}$ be the set of the affine
equilibrium distributions (see \cite{V3}, p.261).
This set spans the kernel of $-(1/2)I +\Cal{K}_\rho$ on
$W^{-1,2}(\partial\OO)\times L^2(\partial\OO)$.
It follows from (7.16) and duality that for $p$ close to $2$,
$$
\aligned
\frac12 I +\Cal{K}^*_\rho&: \ \ 
W^{1,p}(\partial\OO)\times L^p(\partial\OO)
\to W^{1,p}(\partial\OO)\times L^p(\partial\OO),\\
-\frac12 I+\Cal{K}^*_\rho&:\ \ 
\bold{Z}^p(\partial\OO)
\to \bold{Z}^p(\partial\OO),
\endaligned
\tag 7.17
$$
are isomorphisms,
where $\bold{Z}^p(\partial\OO)$ is a subspace of $W^{1,p}(\partial\OO)
\times L^p(\partial\OO)$ whose elements $(F,g)$ satisfy
$$
\Lambda_j^*(F) +\int_{\partial\OO} f_j^*\, g\, d\sigma
=0\ \ \ \text{ for }\ \ j=0,1,\dots, n.
\tag 7.18
$$
Note that $\bold{Z}^p(\partial\OO)$ is well defined for $p>2-\e$.

\proclaim{\bf Theorem 7.1}
There exists $\e>0$ such that
the operators in (7.17) are isomorphisms
for $2<p<\frac{2(n-1)}{n-3} +\e$ and $n\ge 4$.
If $n=2$ or $3$, the operators in (7.17) are isomorphisms
for any $2<p<\infty$.
\endproclaim

Theorem 7.1 follows from Theorem 6.7 
by the same line of argument that we used to
prove Theorem 3.1.
To carry out the proof, we need to compute the Neumann
trace of the double layer potential.
Let $WA_2^p(\partial\OO)$ denote the space of Whitney
arrays $\dot{f}=\left\{ f_0, f_1,\dots, f_n\right\}
\subset W^{1,p}(\partial\OO)$ which satisfy
the compatibility conditions
$\frac{\partial f_0}{\partial T_{ij}}
=N_i f_j -N_j f_i$ for $1\le i<j\le n$
\cite{V2}.

\proclaim{\bf Lemma 7.2}
Let $\dot{f}=\{ f_0, f_1, \dots, f_n\}\in WA^p_2(\partial\OO)$.
Let $w(x)=\Cal{D}_\rho (F,g)$ with $F=f_0$ and $g=-N_i f_i$. Then
$(\nabla \nabla w)^*\in L^p(\partial\OO)$ and
$$
\big( K_\rho(w)_+, M_\rho(w)_+\big)
 =\big(K_\rho (w)_-, M_\rho(w)_-\big),
\tag 7.19
$$
on $\partial\OO$.
\endproclaim

\demo{Proof} 
Using (7.4) and the compatibility conditions, we have
$$
\aligned
&D_\ell w(x)=\\
&-\int_{\partial\OO}
\left\{ D_i\Gamma^x \cdot \frac{\partial F}{\partial T_{\ell i}}
+D_\ell \Gamma^x\cdot g
+(1-\rho)
\frac{\partial}{\partial T_{ik}}
D_k D_\ell B^x\cdot f_i\right\}\, d\sigma\\
&=\int_{\partial\OO}
\frac{\partial\Gamma^x}{\partial N}\, f_\ell
d\sigma
+\int_{\partial\OO}
\left\{ \Gamma^x \cdot \frac{\partial f_i}{\partial T_{\ell i}}
+(1-\rho) D_kD_\ell B^x
\cdot \frac{\partial f_i}{\partial T_{ik}}\right\}\, 
d\sigma .
\endaligned
\tag 7.20
$$
It follows that
$$
\aligned
&D_jD_\ell w(x)=\\
&\int_{\partial\OO}
\left\{ D_i\Gamma^x\cdot \frac{\partial f_\ell}{\partial T_{ij}}
+D_j \Gamma^x \cdot \frac{\partial f_i}{\partial T_{i\ell}}
+(1-\rho)
D_jD_kD_\ell B^x
\cdot \frac{\partial f_i}{\partial T_{ki}}\right\}\, d\sigma.
\endaligned
\tag 7.21
$$
By \cite{CMM}, this implies $\| (\nabla\nabla w)^*\|_p
\le C\, \sum_i \| \nabla_t f_i\|_p<\infty$.
Also it follows from  (7.6) that
$$
D_jD_\ell w_+-D_jD_\ell w_- =N_i\frac{\partial f_\ell}{\partial T_{ij}}
+N_j \frac{\partial f_i}{\partial T_{i\ell}}
+(1-\rho)N_jN_kN_\ell \frac{\partial f_i}{\partial T_{ki}}.
\tag 7.22
$$
This yields that $M_\rho(w)_+=M_\rho(w)_-$ on $\partial\OO$ by a simple
computation.
To find $K_\rho(w)_\pm
=\frac{\partial}{\partial N}
\Delta w_\pm +(1-\rho)\frac{\partial}{\partial T_{ij}}
\big( N_\ell N_i D_jD_\ell w\big)_\pm$ on $\partial\OO$, we note that
by (7.21),
$$
\Delta w(x)=(1-\rho)\int_{\partial\OO}
D_j \Gamma^x \cdot \frac{\partial f_i}{\partial T_{ji}}\, 
d\sigma.
\tag 7.23
$$
Thus we may write
$$
\frac{\partial\Delta w}{\partial N}
=(1-\rho)\frac{\partial }{\partial T_{\ell j}}
\int_{\partial\OO}
D_\ell \Gamma^x \cdot \frac{\partial f_i}{\partial T_{ji}}\, d\sigma.
\tag 7.24
$$
It then follows from (7.24), (7.21) and (7.6) that
$$
\aligned
&\left[ K_\rho(w)_+ -K_\rho(w)_-\right](\phi)\\
&=(1-\rho)
\int_{\partial\OO}
N_\ell \frac{\partial f_i}{\partial T_{ji}}
\cdot \frac{\partial\phi}{\partial T_{j\ell}}\, d\sigma\\
&\ \ \
+(1-\rho)
\int_{\partial\OO}
N_\ell N_i \left\{ N_m\frac{\partial f_\ell}{\partial T_{mj}}
+N_j \frac{\partial f_m}{\partial T_{m\ell}}
+(1-\rho)
N_jN_kN_\ell \frac{\partial f_m}{\partial T_{km}}
\right\}\,
\frac{\partial \phi}{\partial T_{ji}}\, d\sigma\\
&=(1-\rho)\int_{\partial\OO}
\left\{
N_\ell \frac{\partial f_i}{\partial T_{ji}}\cdot 
\frac{\partial\phi}{\partial T_{j\ell}}
+N_\ell N_i N_m \frac{\partial f_\ell}{\partial T_{mj}}
\cdot \frac{\partial \phi}{\partial T_{ji}}\right\}\, d\sigma\\
&=\int_{\partial\OO}
\left\{ N_iN_j \frac{\partial f_\ell}{\partial T_{j\ell}}
-\frac{\partial f_\ell}{\partial T_{i\ell}}
-N_\ell N_m \frac{\partial f_\ell}{\partial T_{mi}}
\right\}\, D_i \phi\, d\sigma\\
&=0,
\endaligned
$$
where we have used the compatibility condition
$$
N_i \frac{\partial f_\ell}{\partial T_{jk}}
=N_k \frac{\partial f_\ell}{\partial T_{ji}}
-N_j \frac{\partial f_\ell}{\partial T_{ki}}
$$
for $k=\ell$ in the last step.
This finishes the proof.
\enddemo

\demo{\bf Proof of Theorem 7.1}
We will give the proof of the invertibility of $-(1/2)I+\Cal{K}^*_\rho$
on $\bold{Z}^p(\partial\OO)$. The case for $(1/2)I +\Cal{K}_\rho^*
$ on $W^{1,p}(\partial\OO)\times L^p(\partial\OO)$ is 
similar.

Let $(G, h)\in \bold{Z}^p(\partial\OO)$ for some $2<p<\infty$.
Since $-(1/2)I +\Cal{K}_\rho^*$
is invertible on $\bold{Z}^2(\partial\OO)$, there exists
$(F,g)\in \bold{Z}^2(\partial\OO)$ so that
$\big(-(1/2)I+\Cal{K}_\rho^*\big)(F,g)=(G,h)$. Let $u(x)=\Cal{D}
_\rho(F,g)$
be the double layer potential. We will show that
if $n\ge 4$ and $2<p<p_n+\e$, or if $n=2$, $3$ and $2<p<\infty$,
$$
\aligned
\bigg\{ \frac{1}{s^{n-1}}
\int_{B(P,s)\cap \partial\OO}
& |(\nabla u)^*|^p\, d\sigma\bigg\}^{1/p}
\le 
C
\left\{ \frac{1}{s^{n-1}}
\int_{B(P,Cs)\cap \partial\OO}
 |(\nabla u)^*|^2\, d\sigma\right\}^{1/2}\\
&\ \ \ \ \
+C\, \left\{ \frac{1}{s^{n-1}}
\int_{B(P,Cs)\cap \partial\OO}
\big( |\nabla_t G| +|h|\big)^p\, d\sigma\right\}^{1/p},
\endaligned
\tag 7.25
$$
for any $P\in \partial\OO$ and $s>0$ small.
Since $(F,g)=(u_+- u_-, -\frac{\partial u_+}{\partial N}
+\frac{\partial u_-}{\partial N})$,
by covering $\partial\OO$ with a finite number of small
balls,
we obtain
$$
\aligned
\| \nabla_t F\|_p
+\| g\|_p
&\le C\, \| (\nabla u)^*\|_p
\le C\, \big\{
\| (\nabla u)^*\|_2
+\| \nabla_t G\|_p
+\| h\|_p\big\}\\
&\le C\, \big\{
\|\nabla_t F\|_2 
+\| g\|_2+
\|\nabla_t G\|_p +\| h\|_p\big\}\\
&\le C\, \big\{
\| \nabla_t G\|_p
+\| h\|_p\big\}.
\endaligned
\tag 7.26
$$
This shows that $-(1/2)I +\Cal{K}^*_\rho$ is
invertible on $\bold{Z}^p(\partial\OO)$.
Note that by a density argument, we may assume that
$(G,h)=(f_0, -f_iN_i)$ for some $\{ f_0, f_1,\dots, f_n\}
\in WA^2_2(\partial\OO)$. This would imply that
$(F,g)=(\widetilde{f}_0, -\widetilde{f}_i N_i)$
for some $\{ \widetilde{f}_0, \widetilde{f}_1,
\dots, \widetilde{f}_n\}
\in WA^2_2(\partial\OO)$ by \cite{V3} (p.265).
Consequently $(\nabla\nabla u)^*\in L^2(\partial\OO)$
by Lemma 7.2. 

To establish estimate (7.25), we may assume that $P=0$ and
$B(0,r_0)\cap\OO$ is given by (2.2).
Let $Q_0=I_s$ be a surface cube defined in (2.3).
For any subcube  $Q$ of $Q_0$, we choose a function
$\varphi=\varphi_Q\in C_0^2(\br^n)$ such that
$0\le \varphi\le 1$, $\varphi=1$ in $100Q$,
$\varphi=0$ outside of $200Q$, and $|\nabla\varphi|
\le C/r$, $|\nabla\nabla\varphi|\le C/r^2$ where
$r$ is the diameter of $Q$. 
Let
$$
\beta =\frac{1}{|200Q|}
\int_{200Q} G\, d\sigma.
\tag 7.27
$$
Since
$$
W^{1,2}(\partial\OO)\times L^2(\partial\OO)
=\bold{Z}^2(\partial\OO)
\oplus \text{span}\big\{ (1,0), (x_j, -N_j), j=1,\dots, n\big\},
\tag 7.28
$$
there exists $(F_Q, g_Q)\in \bold{Z}^2(\partial\OO)$ and $(\alpha_0,\alpha_1,
\dots, \alpha_n)\in \br^{n+1}$ such that
$$
\aligned
\big( &(G-\beta)\varphi,
-h\varphi -(G-\beta)\frac{\partial\varphi}{\partial N}\big)\\
&\ \ \ =(-\frac12 I +\Cal{K}_\rho^*)
(F_Q, g_Q) +\alpha_0 (1,0)
+ \alpha_j (x_j, -N_j),\\
\| (G-\beta)\varphi&\|_{W^{1,2}(\partial\OO)}
+\| h\varphi + (G-\beta)
\frac{\partial\varphi}{\partial N}\|_2\\
&\sim \| F_Q\|_{W^{1,2}(\partial\OO)}
+\| g_Q\|_2 +\sum_{j=0}^n |\alpha_j|.
\endaligned
\tag 7.29
$$
Let
$
 v(x)=\Cal{D}_\rho (F_Q, g_Q) +\alpha_0 +
\alpha_j x_j
$
and $w=u-v-\beta=\Cal{D}_\rho (F-F_Q, g-g_Q)-\beta$. Note that
$$
(w_-,-\frac{\partial w_-}{\partial N})
=\big( (G-\beta)(1-\varphi),
-h(1-\varphi) +(G-\beta)\frac{\partial \varphi}{\partial N}\big).
\tag 7.30
$$
Thus $w_- =|\nabla w_-|=0$ on $100Q$. Since
$(-(1/2)I+\Cal{K}_\rho^*)(F_Q, g_Q)$ is given by an
array in $WA^2_2(\partial\OO)$, we may deduce that
$(F_Q,g_Q)$ is also given by an array in $WA^2_2(\partial\OO)$.
It follows from Lemma 7.2
that $(\nabla\nabla w)^*\in L^2(\partial\OO)$ and
$\big( M_\rho(w)_+, K_\rho(w)_+\big)
=\big( M_\rho(w)_-, K_\rho(w)_-\big)$ on $\partial\OO$.
This allows us to apply Theorem 6.7. We obtain
$$
\left\{\frac{1}{|Q^\prime|}
\int_{Q^\prime} |(\nabla w)^*|^{p_n}\, d\sigma\right\}
^{1/p_n}
\le C\, 
\left\{\frac{1}{|32Q^\prime|}
\int_{32 Q^\prime}
|(\nabla w)^*|^2\, d\sigma\right\}^{1/2}
\tag 7.31
$$
for any subcube $Q^\prime$ of $Q$.
Since the reverse H\"older inequality
(7.31) is self-improving \cite{Gi}, in the case $n\ge 4$,
this means that there exists
$\e>0$ depending only on $\|\psi\|_\infty$, $n$ and the
constant $C$ in (7.31) so that
$$
\aligned
&\left\{\frac{1}{|Q|}
\int_{Q} |(\nabla w)^*|^{\bar{p}}\, d\sigma\right\}
^{1/\bar{p}}
\le C\,
\left\{\frac{1}{|64Q|}
\int_{64 Q}
|(\nabla w)^*|^2\, d\sigma\right\}^{1/2}\\
&\le C\, \left\{
\frac{1}{|64Q|}
\int_{64 Q}
|(\nabla u)^*|^2\, d\sigma\right\}^{1/2}
+C\,
\left\{
\frac{1}{|64Q|}
\int_{64 Q}
|(\nabla v)^*|^2\, d\sigma\right\}^{1/2},
\endaligned
\tag 7.32
$$
where $\bar{p}=p_n+\e$.

Finally we note that
by (7.29)
$$
\aligned
\left\{
\int_{\partial\OO}
|(\nabla v)^*|^2\, d\sigma\right\}^{1/2}
&\le C\, \left\{
\int_{\partial \OO}
\big(|\nabla_t F_Q|+|g_Q|\big)^2\, d\sigma\right\}^{1/2}
+\sum_{j=1}^n |\alpha_j| \\
&\le 
C\, \left\{\int_{200Q}
\big( |\nabla_t G|+|h|\big)^2\, d\sigma\right\}^{1/2},
\endaligned
\tag 7.33
$$
where we also used the Poincar\'e inequality.
With (7.33) and (7.32),
estimate (7.25) follows by
Theorem 3.2. This completes the proof of Theorem 7.1.
\enddemo

\remark{\bf Remark 7.3}
The $L^p$ Dirichlet problem for the biharmonic equation 
$$
\left\{
\aligned
&\Delta^2 u = 0 \ \ \ \text{ in }\ \ \OO,\\
&u=F\in W^{1,p}(\partial\OO),\ \ \
\frac{\partial u}{\partial N}=g\in L^p(\partial\OO)\ \ 
\text{ on } \ \partial\OO,\\
&(\nabla u)^*\in L^p(\partial\OO),
\endaligned
\right.
\tag 7.34
$$
is uniquely solvable if
$$
\alignedat3
&n=2, 3,& \quad\quad & 2-\e<p\le \infty,\\
&n=4, & \quad \quad & 2-\e<p <6+\e,\\
&n=5,6,7, & \quad\quad & 2-\e<p<4+\e,\\
&n\ge 8, &\quad\quad & 2-\e<p<2+\frac{4}{n-\lambda_n}
+\e, 
\endalignedat
\tag 7.35
$$
where $ \lambda_n=(n+10 +2\sqrt{2(n^2-n+2)})/7$.
See \cite{DKV1, PV1, S3, S5}. The ranges of $p$'s in (7.35)
are known to be sharp in the case $2\le n\le 7$ \cite{PV1}.
This implies that the ranges of $p$'s in Theorem 7.1
are sharp for $n=2,3,4,5$.
\endremark

\proclaim{\bf Corollary 7.4}
Let $2<p<p_n +\e$ for $n\ge 4$ and $2<p<\infty$
for $n=2$ or $3$. The unique solution to the Dirichlet problem (7.34)
for the biharmonic equation with boundary data $(F,g)$
is given by
$$
u(x)=\Cal{D}_\rho \left((\frac12 I+\Cal{K}^*_\rho)^{-1}
(F,g)\right).
\tag 7.36
$$
\endproclaim

By duality and an argument similar to that in the proof
of Theorem 5.1, we may deduce the following from Theorem 7.1.

\proclaim{\bf Theorem 7.5}
There exists $\e>0$ such that the operators $\pm (1/2)I +\Cal{K}_\rho
$ in (7.16) are isomorphism for $n\ge 4$ and $
\frac{2(n-1)}{n+1}-\e
<p<2$. If $n=2$ or $3$, the operators are isomorphism
for $1<p<2$.
\endproclaim

\demo{\bf Proof of Theorem 1.3}
The existence follows from the invertibility of
$-(1/2)I+\Cal{K}_\rho$ on $\bold{X}^p(\partial\OO)$,
while the uniqueness was proved in \cite{V3}, p.273
by constructing a Neumann function.
\enddemo

\bigskip

\centerline{\bf 8. The Classical Layer Potentials on Weighted Spaces}

In this section we consider the classical layer potentials
for Laplace's equation
$\Delta u=0$ in $\OO$. In order to be consistant with our notation
for elliptic systems, we shall use the fundamental
solution for $\Cal{L}=-\Delta$ in the definitions of single and double
layer potentials. It is well known that
the operators $(1/2)I +\Cal{K}: L_0^p(\partial\OO) \to L_0^p(\partial\OO)$ and
$-(1/2)I+\Cal{K}: L^p(\partial\OO)
\to L^p(\partial\OO)$ are isomorphisms for $n\ge 2$ and
$1<p<2+\e$. The case $p=2$ was proved in \cite{V1}, using Rellich
identities as we indicated in Section 1. The sharp
range $1<p<2+\e$ was obtained in \cite{DK1}. This was done by
establishing $L^1$ estimates for solutions of the Neumann and
 regularity problems
with boundary data in the atomic Hardy Spaces.
It follows by duality that $(1/2)I +\Cal{K}^*$ and $-(1/2)I+\Cal{K}^*$
are isomorphisms on $L^p(\partial\OO)/\{ h_0\}$ and $L^p(\partial\OO)$
respectively, where $2-\e_1<p<\infty$ and
$h_0$ is a function which spans the kernel
of $(1/2)I+\Cal{K}$ on $L^2(\partial\OO)$.

With the method in previous sections,
it is possible to recover the sharp $L^p$ invertibility
 in \cite{DK1} without the use
of the Hardy spaces. To do this, we will prove directly that
$(1/2)I+\Cal{K}^*:L^p(\partial\OO)/\{ h_0\}
\to L^p(\partial\OO)/\{ h_0\}$ and $-(1/2)I +\Cal{K}^*:
L^p(\partial\OO)\to L^p(\partial\OO)$ are invertible
for $2-\e_1<p<\infty$.
In fact we shall prove a
stronger result.
Let $\Cal{X}^2(\partial\OO, \omega d\sigma)$ denote the space
of functions $f$ in $L^2(\partial\OO, \omega d\sigma)$ such that
$\int_{\partial\OO} f h_0\, d\sigma =0$.

\proclaim{\bf Theorem 8.1}
Let $\OO$ be a bounded Lipschitz domain in $\br^n$, $n\ge 3$ with
connected boundary. Then there exists $\delta\in (0,1]$ depending only
on $n$ and the Lipschitz character of $\OO$ such that
the operators
$$
\aligned
(1/2) I +\Cal{K}^*&:\ 
\Cal{X}^2(\partial\OO,\omega d\sigma)
\to \Cal{X}^2(\partial\OO, \omega d\sigma),\\
-(1/2) I +\Cal{K}^*&:\ 
L^2(\partial\OO,\omega d\sigma)
\to L^2(\partial\OO, \omega d\sigma),
\endaligned
\tag 8.1
$$
are isomorphisms for any $A_{1+\delta}$ weight $\omega$
on $\partial\OO$.
\endproclaim

We refer the reader to \cite{St2} for the theory of $A_p$ weights.
In particular the boundedness of operator $\Cal{K}^*$ on
$L^2(\partial\OO, \omega d\sigma)$ with $\omega\in A_2(\partial\OO)$
follows from \cite{CMM} and the standard weighted inequalities
for Calder\'on-Zygmund
operators.
Also, by H\"older inequality, $L^2(\partial\OO, \omega d\sigma)
\subset L^p(\partial\OO)$ if $\omega\in A_{1+\delta}(\partial\OO)$
 and $p=2/(1+\delta)$. 
Since $h_0\in L^q(\partial\OO)$ for some $q>2$,
this implies that
the space $\Cal{X}^2(\partial\OO, \omega d\sigma )$ is well defined if
$\omega\in A_{1+\delta}$ and $\delta>0$ is sufficiently small.

Note that
by an extrapolation theorem of Rubio de Francia (see e.g. \cite{Du}),
 Theorem 8.1 yields the 
$L^p$ inveribility of $\pm (1/2) I +\Cal{K}^*$
for the sharp range $2-\e<p<\infty$.
Furthermore,
by duality, we obtain the following.

\proclaim{\bf Theorem 8.2}
Let $\OO$ be a bounded Lipschitz domain in $\br^n$, $n\ge 3$ with
connected boundary. Then there exists $\delta\in (0,1]$ depending only
on $n$ and the Lipschitz character of $\OO$ such that
the operators $(1/2)I +\Cal{K}$ and $-(1/2)I+\Cal{K}$
are isomorphisms on $L_0^2\left(\partial\OO,\frac{ d\sigma}{\omega}\right)$
and $L^2\left(\partial\OO, \frac{d\sigma}{\omega}\right)$
respectively, for any $A_{1+\delta}$ weight $\omega$
on $\partial\OO$.
\endproclaim

Here $L^2_0\left(\partial\OO, \frac{ d\sigma}{\omega}\right)$
 denotes the space
of functions $f$ in $L^2\left(\partial\OO, \frac{d\sigma}{\omega}\right)$
 such that
$\int_{\partial\OO} f\, d\sigma =0$. 
To prove Theorem 8.2, one uses the fact that
$L^2(\partial\OO, \omega d\sigma)=\Cal{X}^2(\partial\OO, \omega d\sigma)
\oplus \br$ and preceeds as in the proof of Theorem 4.1.

As in the $L^p$ case, the invertibility of $(1/2)I+\Cal{K}$
on $L^2\left(\partial\OO, \frac{d\sigma}{\omega}\right)$
gives us the existence for
 the Neumann problem with boundary data in the
weighted $L^2$ space. 
Since  $L^2\left(\partial\OO, \frac{d\sigma}{\omega}\right)
\subset L^p(\partial\OO)$ for some $p>1$.
The uniqueness follows from the uniqueness for the 
$L^p$ Neumann problem \cite{DK1}.

\proclaim{\bf Corollary 8.3}
Let $\OO$ be a bounded Lipschitz domain in $\br^n$, $n\ge 3$ with
connected boundary. Then there exists $\delta\in (0,1]$ depending only
on $n$ and the Lipschitz character of $\OO$ such that
given any $g\in L_0^2\left(\partial\OO, \frac{d\sigma}{\omega}\right)$
with $\omega\in A_{1+\delta}(\partial\OO)$, there exists
a harmonic function $u$ on $\OO$, unique up to constants,
such that $\frac{\partial u}{\partial N} =g$ and
$(\nabla u)^*\in L^2\left(\partial\OO,\frac{d\sigma}{\omega}\right)$.
Moreover, the solution $u$ satisfies 
$$
\| (\nabla u)^*\|_{L^2\left(\partial\OO,\frac{d\sigma}{\omega}\right)}
\le C\, \| g\|_{L^2\left(\partial\OO,\frac{d\sigma}{\omega}\right)},
\tag 8.2
$$
and is given by the single layer potential with density
$((1/2)I +\Cal{K})^{-1}(g)$.
\endproclaim

\remark{\bf Remark 8.4}
The condition $\omega\in A_{1+\delta}$ in Theorems 8.1 and 8.2
(and in Corollary 8.3)
is sharp in the context of $A_p$ weights.
 This is because they imply the
sharp ranges of $p$'s for the $L^p$ invertibility.
However in the case $n\ge 4$,
there are weights $\omega$ which are not in the 
sharp $A_p$ class and for which $\pm (1/2)I+\Cal{K}$ are 
invertible
on $L^2\left(\partial\OO, \frac{d\sigma}{\omega}\right)$.
Indeed, consider the power weight $\omega_\alpha
=|Q-Q_0|^\alpha$, where $Q_0\in \partial\OO$ and $\alpha>1-n$.
It is shown in \cite{S2} that $(1/2)I +\Cal{K}$ and
$-(1/2)I +\Cal{K}$ are invertible on $L_0^2\left(\partial\OO, 
\frac{d\sigma}{\omega_\alpha}\right)$ and
$L^2\left(\partial\OO, 
\frac{d\sigma}{\omega_\alpha}\right)$ 
respectively, if $1-n<\alpha<n-3 +\e$.
However we observe that $\omega_\alpha\in A_{1+\delta}$ if and only if
$1-n<\alpha<(n-1)\delta$.
\endremark

It remains to prove Theorem 8.1.
To do this, we need to establish a reverse H\"older
inequality similiar to (2.28), but with $p_n$ replaced
by any exponent $p>2$.
Since $|\nabla u|$ on the boundary is only $L^q$
integrable for some $q>2$, the Sobolev inequality
is not useful in higher dimensions. Instead we use the following
Morrey space estimate (see e.g. \cite{Gi}, Ch.3),
$$
\aligned
&\sup_{I(P_0,R)} |u|
\le 
\\
&\frac{C}{R^{n-1}}
\int_{I(P_0,2R)}
|u|\, d\sigma
+C_\lambda \, R^{\frac{\lambda -n+3}{2}}
\sup\Sb 0<r<R\\
P\in I(P_0,R)\endSb
\left\{ r^{-\lambda}
\int_{I(P,r)}
|\nabla_t u|^2\, d\sigma\right\}^{1/2}
\endaligned
\tag 8.3
$$
where $\lambda>n-3$ and $I(P,r)=B(P,r)\cap\partial\OO$ for
$P\in \partial\OO$ and $0<r<r_0$.

Assume $0\in \partial\OO$ and $\OO\cap B(0,r_0)$ is given by (2.2).

\proclaim{\bf Lemma 8.5}
Suppose $\Delta u=0$ in $\OO_\pm$. Assume that
$(\nabla u)^*_\pm \in L^2(I_{4R})$ and $u_\pm =0$ on $I_{4R}$
for some $0<4R<cr_0$.
Then there exists $\lambda>n-3$ depending only on $n$ and $\OO$ such that
$$
\sup_{0<r<R}
r^{-\lambda}
\int_{I_r} |\nabla u_\pm |^2\, d\sigma
\le \frac{C}{R^{\lambda +3}}
\int_{D^\pm_{4R}}| u|^2\, dx.
\tag 8.4
$$
\endproclaim

\demo{Proof} Since $u_\pm =0$ in $I_{4R}$,
we may use (2.10) and (2.5) to obtain
$$
\int_{I_r} |\nabla u_\pm |^2\, d\sigma
\le \frac{C}{r^3}
\int_{D_{4r}^\pm } |u|^2\, dx.
\tag 8.5
$$
By the boundary H\"older estimates, we have
$$
|u(x)|^2\le C \left(\frac{r}{R}\right)^\delta \frac{1}{R^n}
\int_{D_{4R}^\pm}
|u|^2\, dx,
\tag 8.6
$$
for any $x\in D_{4r}^\pm $, where $\delta>0$ depends only
on $n$ and $\OO$. Estimate (8.4) with $\lambda =n-3+\delta$
now follows easily from (8.5) and (8.6).
\enddemo  

\proclaim{\bf Lemma 8.6}
Suppose that $\Delta u=0$ in $\OO_\pm$ and $(\nabla u)_\pm^*\in L^2(I_{4R})$
for some $0<4R<cr_0$. Then there exists $\lambda>n-3$
depending only on $n$ and $\OO$ such that
$$
\aligned
&\sup_{0<r<R}
r^{-\lambda} \int_{I_r}|\nabla u_\pm |^2\, d\sigma\\
& \le C\, \sup_{0<r<2R}
r^{-\lambda}
 \int_{I_r} |\frac{\partial u_\pm}{\partial N}|^2\,
d\sigma
 +\frac{C}{R^{\lambda +3}}
\int_{D_{4R}^\pm}
|u|^2\, dx
+\frac{C}{R^{\lambda +1}}\int_{I_{4R}}
|u_\pm |\, \big|\frac{\partial u_\pm}{\partial N}\big|\, 
d\sigma. 
\endaligned
\tag 8.7
$$
\endproclaim

\demo{Proof}
We use the following estimate established in \cite{S2} (Lemma 4.18, p.2855),
$$
\int_{I_r} |\nabla u |^2\, d\sigma
\le C\, r^{\lambda_0}
\int_{\partial D_{sR}^\pm } \frac{\big|\frac{\partial u}
{\partial N}\big|^2}
{\{ |P| + r\}^{\lambda_0}}\, d\sigma(P),
\tag 8.8
$$
where $n-3<\lambda_0<n-3+\e$. It follows that
if $n-3<\lambda<\lambda_0$,
$$
\sup_{0<r<R}
r^{-\lambda} \int_{I_r}
|\nabla u|^2\, d\sigma
\le C\sup_{0<r<2R}
r^{-\lambda} \int_{I_r} \big|\frac{\partial u}
{\partial N}\big|^2\, d\sigma
+\frac{C}{R^{\lambda}} \int_{\OO_\pm \cap \partial D_{sR}^\pm}
|\nabla u|^2\, d\sigma,
\tag 8.9
$$
for $1<s<2$.
Estimate (8.7) now follows by an integation in $s$
over $(1,2)$ and using (2.5).
\enddemo

\proclaim{\bf Lemma 8.7}
Suppose that $\Delta u=0$ in $\br^n\setminus \partial\OO$
and $(\nabla u)^*_+ +(\nabla u)^*_-\in L^2(I_{16R})$
for some $0<16R<cr_0$.
Assume that either $u_+=0$ or $u_-=0$ on $I_{16R}$.
Then there exists $\lambda>n-3$ and $p_0<2$ depending only
on $n$ and $\OO$ such that
$$
\aligned
\sup_{0<r<R} r^{-\lambda} \int_{I_r} |\nabla u_\pm|^2\, d\sigma
&\le C\, \sup_{0<r<2R} r^{-\lambda}\int_{I_r}
\big| \frac{\partial u_+}{\partial N} 
-\frac{\partial u_-}{\partial N}\big|^2\, d\sigma\\
&\ \ \ \  +\frac{C}{R^{\lambda +1}}
\int_{I_{8R}}
\big( |u_+| +|u_-|\big) \big| \frac{\partial u_+}{\partial N}
-\frac{\partial u_-}{\partial N}\big|\,
d\sigma\\
&\ \ \ \ +\frac{C}{R^{\lambda +3}}\int_{D_{16R}^+\cup 
D_{16R}^-} |u|^2\, dx\\
&\ \ \ \ +C\, R^{n-\lambda-3}
\left\{ \frac{1}{R^{n-1}}
\int_{I_{8R}}
\big( |u_+| +|u_-|\big)^{p_0}\, d\sigma\right\}^{2/p_0}.
\endaligned
\tag 8.10
$$
\endproclaim

\demo{Proof}
We only consider the case $u_+=0$ on $I_{16R}$.
The case for $u_-$ is exactly the same.
 
The estimate for $r^{-\lambda} \int_{I_r} |\nabla u_+|^2 d\sigma$
is contained in (8.4). To estimate $r^{-\lambda} \int_{I_r}
|\nabla u_-|^2d\sigma$, in view of (8.7) and (8.4), we only need
to take care of the term
$$
\frac{1}{R^{\lambda +1}}
\int_{I_{4R}} |u_-|\, \big|\frac{\partial u_-}
{\partial N}\big|\, d\sigma.
\tag 8.11
$$
To this end,
first we  replace $\big|\frac{\partial u_-}{\partial N}\big|$ in (8.11) by
$\big|\frac{\partial u_+}{\partial N}\big|$, since the 
difference is bounded by the second term in the right side of (8.10).
Next we use the H\"older inequality. This reduces the problem to the
estimation of
$$
R^{n-\lambda-1}\left\{ \frac{1}{R^{n-1}}
\int_{I_{4R}} \big|\frac{\partial u_+}{\partial N}\big|^{p_0^\prime}
d\sigma\right\}^{2/p_0^\prime}.
\tag 8.12
$$
Finally we use the $L^{p_0^\prime}$ estimate for the regularity problem
on $D_{sR}^+$ for $s\in (4,5)$ and then a familiar integation in $s$
to bound the term in (8.12) by
$$
\aligned
C\, R^{n-\lambda-1}\left\{\frac{1}{R^n}
\int_{D_{5R}^+}
|\nabla u|^{p_0^\prime} dx\right\}^{2/p_0^\prime}
& \le \frac{C}{R^{\lambda +1}}
\int_{D_{6R}^+}
|\nabla u|^2\, dx\\
&\le \frac{C}{R^{\lambda +3}}
\int_{D_{16R}^+}
|u|^2\, dx,
\endaligned
\tag 8.13
$$
where we have used a higher integrability estimate
in the first inequality (see e.g. \cite{Gi}).
We remark that $L^{p_0^\prime}$ regularity
estimate holds if $p_0$ is close to $2$ \cite{DK1}.
This completes the proof of (8.10).
\enddemo

We now are ready to prove the desired reverse H\"older
inequality.

\proclaim{\bf Theorem 8.8}
Suppose that $\Delta u=0$ in $\br^n\setminus \partial\OO$ and
$(\nabla u)^*_+ +(\nabla u)_-^* \in L^2(I_{300R})$ for some $0<300R<cr_0$.
Also assume that
$\frac{\partial u_+}{\partial N}=
\frac{\partial u_-}{\partial N}$ on $I_{300R}$
and  that either $u_+=0$ or $u_-=0$ on $I_{300R}$.
Then 
for any $2<q<\infty$,
$$
\left\{ \frac{1}{R^{n-1}}
\int_{I_R} |(u)^*|^q \, d\sigma\right\}^{1/q}
\le C_q\, \left\{\frac{1}{R^{n-1}}
\int_{I_{300R}}
|(u)^*|^{p_0}d\sigma\right\}^{1/p_0},
\tag 8.14
$$
where
$p_0<2$ depends only on $n$ and $\OO$.
\endproclaim

\demo{Proof}
It follows from (8.3) and (8.10) that
$$
\aligned
&\left\{\frac{1}{R^{n-1}}
\int_{I_R}
\big( |u_+| +|u_-|\big)^q d\sigma\right\}^{1/q}\\
&\le
C\left\{ \frac{1}{R^n}
\int_{D_{32R}^+\cup D_{32R}^-}
|u|^2 dx\right\}^{1/2}
+C\left\{ \frac{1}{R^{n-1}}
\int_{I_{16R}}
\big( |u_+| +|u_-|\big)^{p_0} d\sigma\right\}^{1/p_0}\\
&\le C\left\{ \frac{1}{R^{n-1}}
\int_{I_{64 R}}
|(u)^*|^{p_0} d\sigma\right\}^{1/p_0},
\endaligned
\tag 8.15
$$
where we also used (2.14) for the second inequality.
Since the $L^p$ Dirichlet problem for Laplace's equation
is solvable for any $p\ge 2$
(this follows from the $L^2$ solvability and the
maximum principle), estimate (8.14)
follows from (8.15) and (2.24).
\enddemo

\demo{\bf Proof of Theorem 8.1}
We only give the proof for the invertibilty of $(1/2)I+\Cal{K}^*$
on $\Cal{X}^2(\partial\OO, \omega d\sigma)$. The case of 
$-(1/2)I +\Cal{K}^*$ is similar.

Let $f\in \Cal{X}^2(\partial\OO, \omega d\sigma)
\cap W^{1,2}(\partial\OO)$. Since $(1/2)I +\Cal{K}^*$
is invertible on $W^{1,2}(\partial\OO)/\{ h_0\}$
and $L^2(\partial\OO)/\{ h_0\}$
\cite{V1}, there exists $g\in W^{1,2}(\partial\OO)$
such that $((1/2)I +\Cal{K}^*) g=f$
and $\| g\|_2\le C\, \| f\|_2$.
We need to show that
$$
\int_{\partial\OO} |g|^2\, \omega d\sigma
\le C\, \int_{\partial\OO} |f|^2\, \omega d\sigma.
\tag 8.16
$$
To this end, we fix $P_0\in \partial\OO$ and $s>0$ sufficiently 
small. Let $u=\Cal{D}(g)$.
We will show that
there exists $p_0<2$ such that
$$
\aligned
\left\{
\int_{I(P_0,s)}
|(u)^*|^2\, \omega d\sigma\right\}^{1/2}
&\le C\, \big\{ \omega(I(P_0, Cs))\big\}^{1/2}
\left\{ \frac{1}{s^{n-1}}\int_{I(P_0,Cs)}
|(u)^*|^{p_0} d\sigma\right\}^{1/p_0}\\
&\ \ \ \ \ \ \ \ \ 
+C\left\{
\int_{I(P_0,Cs)} 
|f|^2 \, \omega d\sigma\right\}^{1/2},
\endaligned
\tag 8.17
$$
for all $\omega\in A_{2/p_0}(\partial\OO)$. 
Note that $\| g\|_{p_0}\le C\, \| f\|_{p_0}$ if 
$p_0$ is close to $2$. Thus
 the first term in the right side of (8.17)
is bounded by
$$
C_s \, \big\{ \omega(\partial\OO)\big\}^{1/2}
\| g\|_{p_0}
\le C_s \, \big\{ \omega(\partial\OO)\big\}^{1/2}
\| f\|_{p_0}
\le C_s \, \| f\|_{L^2(\partial\OO, \omega d\sigma)}.
\tag 8.18
$$
Since $|g|\le 2(u)^*$, estimate (8.16) follows from (8.17) and (8.18)
by covering $\partial\OO$ with a finite number of small surface balls.

We will use Theorem 3.4 to
 prove (8.17).
We may assume that $P_0=0$ and $B(0, r_0)\cap\OO$
is given by (2.2).
Let $Q$ be a small subcube of $I_s$. We proceed as in the proof of Theorem
3.1 to choose function $\varphi =\varphi_Q\in C_0^1 (\br^n)$ and then
$g_Q$ so that $ f\varphi =((1/2)I+\Cal{K}^*)(g_Q) +b$
and $\| f\varphi\|_{p_0} \sim \| g_Q\|_{p_0} +|b|$.
Let
$$
F=|(u)^*|^{p_0},\ \ \ 
R_Q =2^{p_0-1} |(w)^*|^{p_0},
\ \ \text{ and }\ \ 
F_Q=2^{p_0-1} |(v)^*|^{p_0},
\tag 8.19
$$
where $p_0<2$ is given in Theorem 8.8,
 $v=\Cal{D}(g_Q) +b$ and $w=u-v$. Since $w_-=f(1-\varphi)$
and $\frac{\partial w_+}{\partial N}
=\frac{\partial w_-}{\partial N}$, by Theorem 8.8,
we have 
$$
\left\{\frac{1}{|2Q|}\int_{2Q}
|R_Q|^p\, d\sigma\right\}^{1/p}
\le \frac{C}{|Q|}\int_{600Q}| R_Q|\, d\sigma
\tag 8.20
$$
for any $p>(2/p_0)$. Also note that
$$
\| F_Q\|_1 =\| (v)^*\|_{p_0}^{p_0}
\le C\, \big\{ \| g_Q\|_{p_0} +|b|\big\}^{p_0}
\le C\| f\varphi\|_{p_0}^{p_0}.
\tag 8.21
$$
This shows that conditions (3.3) and (3.4) in Theorem 3.2
hold for any $1<p<\infty$. It then follows from Theorem 3.4
and Remark 3.5 with $q=(2/p_0)$ that
estimate (8.17) holds for any $w\in A_{2/p_0}(\partial\OO)$.
This completes the proof.
\enddemo

\remark{\bf Remark 8.5} If $\omega\in A_{1+\delta}(\partial\OO)$,
the Dirichlet problem for Laplace's equation with boundary data
in $L^2(\partial\OO, \omega d\sigma)$
is uniquely solvable. This follows easily
from \cite{D}. In \cite{S2}, we solved the regularity problem
with data in $W^{1,2}\left(\partial\OO, \frac{d\sigma}{\omega}\right)$
for $\omega\in A_{1+\delta}(\partial\OO)$,
and established the sharp estimate
$$
\|(\nabla u)^*\|_{L^2\left(\partial\OO,\frac{d\sigma}{\omega}\right)}
\le C\, \| \nabla_t u\|_{L^2\left(\partial\OO,\frac{d\sigma}{\omega}\right)}.
\tag 8.22
$$
This, together with (8.2), gives the Rellich estimate
(1.24) in the weighted $L^2$ space.
\endremark

\enddocument